\documentclass[11pt]{amsart}
\usepackage[latin1]{inputenc}
\usepackage{amsmath,amsthm}
\usepackage{amssymb,amsfonts}
\usepackage{hyperref}
\usepackage{tikz}

\usepackage[shortlabels]{enumitem}
\usepackage{bm}
\usepackage{cite}
\usepackage{graphicx}
\usepackage{verbatim}
\usepackage{amsthm}
\usepackage{setspace}
\usepackage{amsmath}
\usepackage{color}

\makeatletter
\newcommand{\labeltext}[2]{%
  \@bsphack
  \csname phantomsection\endcsname 
  \def\@currentlabel{#1}{\label{#2}}%
  \@esphack
}
\makeatother

\usepackage[a4paper]{geometry}

\hypersetup{bookmarks,final,colorlinks=true,
linkcolor=[rgb]{0.1, 0.1, 0.44},
citecolor=[rgb]{0.1, 0.1, 0.44}}

\usepackage[mode=multiuser,status=draft,lang=english]{fixme}
\fxsetup{theme=colorsig}

\FXRegisterAuthor{M}{aM}{MV}
\FXRegisterAuthor{G}{aG}{GV}

\usepackage{stmaryrd}
\DeclareSymbolFont{bbold}{U}{bbold}{m}{n}
\DeclareSymbolFontAlphabet{\mathbbm}{bbold}

\title[Model companionship, the continuum problem, and the universe of sets]{What model companionship can say about the continuum problem}

\author{Giorgio Venturi, Matteo Viale}
\thanks{The first author  acknowledges support from the Jovem Pesquisador FAPESP project n. 2016/25891-3 and from the RyC grant  USE-21881-P. The second author acknowledges support from the project:
\emph{PRIN 2017-2017NWTM8R
Mathematical Logic: models, sets, computability} and GNSAGA.
}

	\newtheorem{theorem}{Theorem}[section]
	
	\newtheorem{lemma}[theorem]{Lemma}
	
	\newtheorem{fact}[theorem]{Fact}

\theoremstyle{definition}

	\newtheorem{definition}[theorem]{Definition}
	\newtheorem{notation}[theorem]{Notation}

\theoremstyle{remark}
	\newtheorem{remark}[theorem]{Remark}

\newcommand{\Ord}{\ensuremath{\mathrm{Ord}}}

\newcommand{\ZFC}{\ensuremath{\mathsf{ZFC}}}
\newcommand{\ZF}{\ensuremath{\mathsf{ZF}}}

\DeclareMathOperator{\dom}{dom}
\DeclareMathOperator{\ran}{ran}

\newcommand{\NS}{\ensuremath{\mathbf{NS}}} 

\newcommand{\bool}[1]{\mathsf{#1}}

\newcommand{\SpecAMC}[2]{\mathfrak{spec}_{\mathsf{AMC}}\left(#1,#2\right)}
\newcommand{\SpecMC}[2]{\mathfrak{spec}_{\mathsf{MC}}\left(#1,#2\right)}
\newcommand{\pow}[1]{\mathcal{P}\left(#1\right)}

\newcommand{\qp}[1]{\left[ #1 \right]}
\newcommand{\Qp}[1]{\left\llbracket #1 \right\rrbracket}
\newcommand{\ap}[1]{\langle #1 \rangle}
\newcommand{\bp}[1]{\left\lbrace #1 \right\rbrace}

\newcommand{\BMM}{\ensuremath{\text{{\sf BMM}}}} 
\newcommand{\MM}{\ensuremath{\text{{\sf MM}}}} 
\newcommand{\AX}{\ensuremath{\text{{\sf AX}}}} 
\newcommand{\CH}{\ensuremath{\text{{\sf CH}}}} 

\newcommand{\SSP}{\ensuremath{\text{{\sf SSP}}}}

\newcommand{\A}{{\sf A}}

\begin{document}
\begin{abstract}
We present recent results on the model companions of set theory, placing them in the context of a current debate in the philosophy of mathematics. 
We start by describing the dependence of the notion of model companionship on the signature, and then we analyze this dependence in the specific case of set theory.  We argue that the most natural model companions of set theory describe (as the signature in which we axiomatize set theory varies) theories of $H_{\kappa^+}$, as $\kappa$ ranges among the infinite cardinals. We also single out $2^{\aleph_0}=\aleph_2$ as the unique solution of the continuum problem which can (and does) belong to some model companion of set theory (enriched with large cardinal axioms). 
While doing so we bring to light that set theory enriched by large cardinal axioms in the range of supercompactness has as its model companion (with respect to its first order axiomatization in certain natural signatures) the theory of $H_{\aleph_2}$ as given by a strong form of Woodin's axiom $(*)$ (which holds assuming $\bool{MM}^{++}$).
Finally this model-theoretic approach to set-theoretic validities is explained and justified in terms of a form of maximality inspired by Hilbert's axiom of completeness.  
\end{abstract}

\maketitle


\section{Introduction} \label{sec:introduction}

There is no doubt that the continuum problem is one of the driving forces of set theory. The attempts to determine the cardinality of the continuum  has accompanied the history of set theory, from its very beginning to the present day, and motivated many of its most significant advances. Cantor's definition of the Perfect Set Property was prompted by a partial solution to the continuum Problem; more in general the initial developments of descriptive set theory were also driven by an attempt to confirm the continuum Hypothesis ($\CH$: $2^{\aleph_0}=\aleph_1$) at least for the definable subsets of $\mathbb{R}$. In more recent times, G\"odel's constructible universe and Cohen's method of forcing were devised to show the independence of $\CH$ from $\ZFC$ (the standard first order axiomatization of set theory). 

 The techniques developed by G\"odel and Cohen clarified the intrinsic limitations of the axiomatic approach to the continuum problem, and profoundly influenced the subsequent development of set theory. As a matter of fact, after the 1960's, the many independence results obtained combining the methods of inner models and forcing partially shifted the main focus of set theory: set theorists progressively devoted more and more efforts to the study of the models of set theory, and only derivatively to the study of sets (understood as autonomous mathematical objects).

The present paper presents a new approach to the quest of a solution of the continuum problem and -more generally- to the search of new set theoretic truths; this approach stems from a model-theoretic perspective on set theory, and aims to clarify the philosophical import of the results in \cite{VIAVENMODCOMP,viale2021absolute}. As is standard, we recognize the intrinsic limitations of $\ZFC$ in capturing set-theoretic validities, but, to overcome these limitations,
we bring forward the idea that a powerful criterion to validate a new axiom is given by the nice model theoretic properties of the first order presentations of set theory which include it.
%
 
 The central notion we borrow from model theory and that guides the present approach is that of model companionship. This concept was developed by Abraham Robinson in the 1960's, and it is meant to capture, in an abstract setting, the closure properties that algebraically closed fields display with respect to commutative rings with no zero divisors. 
 As algebraically closed fields contain solutions for all Diophantine equations, the 
 models of the model companions of set theory will have solutions for all ``simple'' set theoretical problems. 
 
 A fundamental property of 
 model companionship is its dependence on the signature.\footnote{See Section \ref{sec:exclosmod} for a precise formulation of this dependence.} Loosely speaking: there can be distinct first order axiomatizations of a mathematical theory $T$, one in signature $\sigma$ and the other in signature $\tau$, such that $T$ admits a model companion when axiomatized according to $\sigma$, but does not when axiomatized according to $\tau$.\footnote{This has an analogy with what occurs for Birkhoff's 
characterization of an algebraic variety in terms of  a universal theory - the latter axiomatization may only be expressed in a signature different from the one in which the variety is initially axiomatized.}
 
The dependence on the signature of model companionship motivates a detailed  discussion aimed to single out the relevant signatures for set theory. We will start by clarifying some fundamental aspects of set-theoretical practice. In particular we will explain why bounded formulae  (i.e. $\Delta_0$-formulae) express set-theoretic concepts of low complexity and therefore should be included in any reasonably rich signature for set theory. We will then show how a careful choice of the signatures for set theory allows one to classify the complexity of set-theoretic concepts and to clarify the informal notion of ``simplicity of a set-theoretic assertion''. Our analysis will provide precise criteria able to link logic complexity,
 forcing invariance, and  simplicity of set-theoretic concepts. Once we have outlined clear criteria to select signatures, we will study set theory under different expansions of $\{\in\}$. By varying the signature,  we will show that the corresponding model companions of set theory describe possible theories of the structures $H_{\kappa^+}$, for $\kappa$ an infinite cardinal.\footnote{Recall that  $H_{\kappa^+}$ is the collection of sets whose transitive closure has size at most  $\kappa$.} Moreover, given an appropriate signature $\tau_\kappa$ for set theory (relative to a cardinal $\kappa$), the corresponding $\tau_{\kappa}$-theory of  $H_{\kappa^+}$ will describe a structure which is closed-off with respect to specific constraints. Roughly speaking: much alike an algebrically closed field $K$ contains solutions to all polynomial equations with coefficients in $K$, the $\tau_\kappa$-structures realizing the model companion of set theory in signature $\tau_\kappa$ will make true all $\Sigma_1$-properties for $\tau_\kappa$ (with the relevant parameters from $H_{\kappa^+}$) which are not ``outward inconsistent''.
 
A precise syntactic description of this closing-off procedure will be given by the notion of absolute model companionship (AMC). We defer to \S \ref{subsec:AMC} a thorough presentation of this notion, we can briefly describe it rightaway by saying that 
AMC describes those $\tau$-theories $S$ for which a model companion exists and is axiomatized by the $\Pi_2$-sentences for $\tau$ which are consistent with the universal and existential fragment of any completion of $S$.\footnote{AMC appears to be the correct notion of model companionship to apply to set theory. For example it is unknown to the authors whether there can be a signature for which the axiomatization of set theory in it admits a model companion which is not its AMC.}
 
 Building on the properties of AMC, the sought solution to the continuum problem can be presented as follows:
\begin{enumerate}[(i)] 
\item  \label{item(i)notCH}
$\neg\CH$ and even a definable version\footnote{Thm. \ref{mainthm:2omegageqomega2*} and the discussion surrounding it will clarify what is meant here by the expression ``a definable version of $2^{\aleph_0}=\aleph_2$''.} $\theta$ of $2^{\aleph_0}=\aleph_2$ can be expressed by $\Pi_2$-formulae in any signature for set theory which contains the $\Delta_0$-properties among its atomic formulae and has a parameter for $\omega_1$.
\item \label{item(ii)notCH}
For each infinite (and definable) cardinal $\kappa$ there is at least one (actually several) signature(s) $\tau_\kappa$ containing the $\Delta_0$-properties among its atomic formulae, a constant to interpret the cardinal $\kappa$, and such that
set theory axiomatized in the corresponding $\tau_\kappa$ admits an AMC.\footnote{Furthermore these AMCs describes theories of $H_{\kappa^+}$ (in ways which we will unfold later on).}
\item  \label{item(iii)notCH}
Given any signature $\tau$ for set theory, containing the $\Delta_0$-properties among its atomic formulae, any  first order $\tau$-axiomatization $T$ of set theory admitting an AMC $S$ for $\tau$ will not have neither $\CH$ nor $\neg\theta$ among the $\tau$-theorems of $S$, unless $T$ models $\CH$ or $\neg\theta$ (furthermore for many choices of $\tau$ and $T$, even if $T$ models $\CH$ and/or $\neg\theta$, $S$ will not).
\item \label{item(iv)notCH}
There is at least one (recursive and natural\footnote{The discussion to follow makes clear that $\tau^*$ is a signature which naturally stems out from the set theoretic practice of the last fourty years.}) signature $\tau^*$ for set theory which contains the $\Delta_0$-properties among its atomic formulae and such that all $\tau^*$-axiomatizations $T$ of \emph{any} set theory enriched with large cardinal axioms\footnote{E.g. $T$ could have $\CH$ or $2^{\aleph_0}>\aleph_2$ among its axioms provided it extends $\ZFC+$\emph{there are class many supercompact cardinals}.} admit an AMC $S$ for $\tau^*$ that contains $\neg\CH$ and $\theta$ among its axioms.
Furthermore $S$ describes the $\tau^*$-theory of $H_{\aleph_2}$ in models of strong forcing axioms. 
\end{enumerate}
 Therefore, these results single out (a definable version of) $2^{\aleph_0}=\aleph_2$ as the unique solution of the continuum problem which can fall in at least one AMC of set theory enriched with large cardinal axioms (for example the AMC with respect to signature $\tau^*$ of item \ref{item(iv)notCH} above).\footnote{Furthermore the same argument we used for $\CH$ and $\theta$ applies mutatis-mutandis to any set theoretic statement $\psi$ which can be formalized by a $\Pi_2$-sentence in the signature $\tau^*$ of item \ref{item(iv)notCH} and such that $\psi$ is forcible over any model of set theory enriched by large cardinal axioms.}

The present model-theoretic approach to set theoretic validities complements well standard strategies for producing axiomatic extensions of $\ZFC$. Actually we see our approach as an implementation of G\"odel's program. Given the stratification of $V$ as the increasing union of its initial segments provided by the various $H_{\kappa^+}$ (as $\kappa$ ranges over the infinite cardinals), we can interpret the AMCs of set theory given by item \ref{item(ii)notCH} as tentative descriptions of the theories of these $H_{\kappa^+}$, theories which coherently extend each other.\footnote{We will come back on this point with more precision (see \ref{cond:Amainres} and the discussion following it, and Thm. \ref{Thm:AMCsettheory+Repl}); roughly the argument is that we can choose uniformly $\tau_\lambda$ for each $\ZFC$-definable cardinal $\lambda$ so that in a model $(V,\in)$ of $\ZFC$ the $\tau_\lambda$-theory of $H_{\lambda^+}^V$  is a model of the AMC of  set theory in signature $\tau_\lambda$ and is naturally described in $H_{\kappa^+}^V$ for $\kappa>\lambda$ by looking at $H_{\lambda^+}$ as a definable subset of $H_{\kappa^+}$.}
Assuming this perspective, these $\tau_\kappa$-theories of $H_{\kappa^+}$ are gradually closing-off the universe of sets with respect to all ``simple'' set theoretical problems (with each $\tau_\kappa$ giving its own description of what is a ``simple'' set theoretic concept, and with each $\tau_\kappa$ making more and more set theoretic concepts becoming ``simple'' as $\kappa$ increases). In this sense, the nice model-theoretic properties displayed by the theories of $H_{\kappa^+}$ (when they describe possible model companions of set theory) realize a form of maximality that we dare to call Hilbertian completeness. As is well-known, Hilbert's axiomatization of geometry expresses a form of maximality for geometrical objects, by stating that the domain of geometry could not be extended if it were to preserve the axioms of the \emph{Grundlagen der Geometrie}. Similarly, the closing-off processes that each model companionship result carries along can (and will) be described as a way to assert the maximization of set-theoretical objects.  
Moreover, important axiomatic extensions of $\ZFC$ (considered candidates for the realization of  G\"odel's program\footnote{A clear justification of why large cardinal axioms are a partial realization of G\"odel's program can be found in \cite{10.1093/philmat/nkj009}. An analysis of why (strong) forcing axioms can also be seen as a realization of this program can be found in \cite{Bag2005} or in the introduction of \cite{VIAMM+++}. See also \cite{WoodinCH1,WoodinCH2,WoodinCHLC} where it is proposed to realize G\"odel's program by introducing new axioms to describe the theories of $H_{\aleph_1}$ and $H_{\aleph_2}$.}) like large cardinals and forcing axioms find, in the present approach, a further justification. As a matter of fact, the generic absoluteness results that we can obtain from large cardinal hypotheses  are now of pivotal importance in detecting the suitable signature expansions able to allow the existence of a(n absolute) model companion for set theory. Another surprising fact which shows that model companionship results and forcing axioms complement each other is that the theory of $H_{\aleph_2}$ in models of strong forcing axioms describes the (absolute) model companion of $\ZFC+$\emph{large cardinals} (as expressed in the natural signature $\tau^*$ of item \ref{item(iv)notCH} above); conversely the existence of a model companion for the theory of $V$  in this signature $\tau^*$ already implies that a strong form of bounded Martin's maximum (which is one of the strongest forcing axioms) and Woodin's axiom $(*)$ both hold in $V$.

\smallskip

The paper is structured as follows:
in \S \ref{sec:rightsignmathth} and \ref{sec:simpcomplconc} we gently analyze the role that signatures can play in outlining the properties of a mathematical theory $T$ (with a focus on how the morphisms between models of $T$ can be used to detect the right signatures for $T$). 
 In \S \ref{sec:exclosmod} we present a few standard notions from model theory (i.e. existentially closed structures, model completeness, and model companionship), which are discussed and generalized in 
 \S \ref{sec:partmorlAMCspec}, where we define the notions of partial Morleyization, absolute model companionship, and (absolute) model companionship spectrum. 
 From \S \ref{sec:AMCspecST} onwards we shift our attention to set theory. We first present precise mathematical criteria to detect what are the ``simple'' set-theoretical concepts. Then
 we analyze various signatures for set theory in \S \ref{subsec:rightsignst}. In \S \ref{subsec:levabsAMCST} we show that the intended models of the model companions of set theory are structures of the form $H_{\kappa^+}$, for $\kappa$ an infinite cardinal.  Finally, in \S \ref{sec:mainresults} we present the main results regarding the model companions of set theory and discuss the information they convey on the cardinality of the continuum. 
We conclude, in \S \ref{sec:philcons}, by comparing our approach with the various notions of maximality we find in the literature and with the debate on the justification of new axioms in set theory. 

\smallskip


This paper aims to reach scholars with interests at the crossroads of philosophy and mathematics. We tried for this reason to minimize the mathematical and philosophical prerequisites needed to follow it. Those who made through this introduction should encounter no serious obstacles to read the remainder of this article. The model companionship results we discuss in the present paper are presented without proofs. The reader interested in them is referred to \cite{VIAVENMODCOMP,viale2021absolute}.

\subsection*{Acknowledgements}
We thank Gabriel Goldberg and the two anonymous referees of RSL for many useful comments on previous drafts of this paper; we also thank David Asper\'o, Ilijas Farah, Leon Horsten, Luca Motto Ros, and Boban Veli\v{c}kovi\'{c} for many interesting discussions on these topics. 
Clearly the opinions expressed (and the errors occurring) in this article are the sole responsibility of the authors and do not extend in any way to those who had the patience to give some advice to improve our presentation.


\section{What is the right signature for a mathematical theory?} \label{sec:rightsignmathth}

One of the great successes of mathematical logic consists in providing an efficient formalization of mathematics: by means of first order logic it is possible to render mathematical theories the objects of a mathematical investigation. In this way logic is able to produce unexpected and non-trivial mathematical results as well as novel insights on a variety of mathematical fields.
It is a matter of facts that there can be many distinct first order  formalizations of a mathematical theory: varying the linguistic presentation of a theory (i.e. its \emph{signature}) we obtain different axiomatic presentations of the same set of theorems. In order to appreciate the variety of possibilities we can encounter, let us consider the concrete case of group theory. This will help us to gently introduce one of the main themes of the present present: the role of the signature in detecting the properties of a theory formalized in it.

We can formalize group theory in first order logic using the signature $\bp{\cdot}$, which simply consists of a binary function symbol and the following axioms:
\begin{align*}
&\forall x,y,z\, [(x\cdot y)\cdot z=x\cdot(y\cdot z)],\\
&\exists x\forall y \,(x\cdot y=y\wedge y\cdot x=y),\\
&\forall x\exists y\forall z\, 
[(x\cdot y)\cdot z=z \wedge (y\cdot x)\cdot z=z  \wedge z\cdot (x\cdot y)=z \wedge z\cdot (y\cdot x)=z].
\end{align*}
Notice that the third axiom, which expresses the existence of a multiplicative inverse, represents a rather complicated assertion, both from the point of view of its syntactic readability and of its L\'evy complexity (being a $\Pi_3$-sentence). The reason is that in this basic signature we lack a constant symbol to denote the neutral element of a group. 
Enriching the language to $\bp{\cdot,e}$, with $e$ a constant symbol, we can now formalize group theory with a simpler set of axioms:
\begin{align*}
&\forall x,y,z\, [(x\cdot y)\cdot z=x\cdot(y\cdot z)],\\
&\forall y \,(e\cdot y=y\wedge y\cdot e=y),\\
&\forall x\exists y\, 
[x\cdot y=e \wedge y\cdot x=e].
\end{align*}
We increased readability and decreased (L\'evy) complexity, as we  are now dealing with a $\Pi_2$-axiomatization. But we can do better. 
Indeed, if we consider the signature $\bp{\cdot,e,{}^{-1}}$, which further adds a unary operation symbol for the inverse operation, we can axiomatize group theory with a 
set of universal equations ($\Pi_1$-sentences):
\begin{align*}
&\forall x,y,z\, [(x\cdot y)\cdot z=x\cdot(y\cdot z)],\\
&\forall x \,(x\cdot e=x\wedge e\cdot x=x),\\
&\forall x\, 
[x\cdot x^{-1}=e \wedge x^{-1}\cdot x=e].
\end{align*}
On the other hand, we could follow a completely different route and axiomatize group theory avoiding the use of function symbols. 
For example we can consider the signature $\bp{R,e}$ consisting of a ternary relation symbol $R$ and a constant symbol $e$ and produce the following axiomatization for group theory:
\begin{align*}
&\forall x,y\exists ! z\, R(x,y,z),\\
&\forall x,y,z,w,t\,[((R(x,y,w)\wedge R(y,z,t))\rightarrow \exists u\, (R(x,t,u)\wedge R(w,z,u))],\\
&\forall y \,[R(e,y,y) \wedge R(y,e,y)],\\
&\forall x\exists y\,[R(x,y,e)\wedge R(y,x,e)]. 
\end{align*}
At the cost of further complicating our axiomatization, we could even drop the use of the constant symbol $e$ and formalize group theory in the signature $\{R\}$. The latter is clearly an artificial solution; moreover, the minimality of the signature does not help the perspicuity of the axiomatization.

Of all the above formalizations, the one mathematicians use more frequently is certainly the one in signature $\bp{\cdot,e,{}^{-1}}$. 
On the contrary, to recognize that the system in signature $\bp{R,e}$ is an axiomatization of group theory would surely require a bit of thought. 

These considerations suggest that among the many possible signatures in which we can formalize a mathematical theory some are better than others. Our aim is to unfold criteria which allow us to detect the best signatures $\tau$ for the formalization of a mathematical theory $T$. To do so we appeal to two sets of arguments.\footnote{(IntCrit) stands for ``Internal criterion'', (AbsCrit) stands for ``Abstract Criterion''.}
\begin{enumerate}[(IntCrit)]
\item \label{INTJUST}
On the one hand, we can select $\tau$ on the basis of specific considerations internal to $T$. For example, by checking the adherence of a $\tau$-formalization of $T$ to the standard informal practice (e.g. the signature $\bp{\cdot, e,{}^{-1}}$  clearly gives the best presentation of group theory in terms of its basic operations). 
\end{enumerate}
\begin{enumerate}[(AbsCrit)]
\item \label{ABSJUST}
On the other hand, we can give abstract criteria for the choice of $\tau$ based on the structural properties that the $\tau$-axiomatization of $T$ displays,  disregarding any consideration on the adherence of the $\tau$-axiomatization of $T$ with the informal one.
\end{enumerate}

It occurs for many $T$ that the reasons internal to $T$ for choosing a signature $\tau$ for its axiomatization cohere with the abstract criteria, and that the two approaches produce the same outcome. This is desirable as abstract criteria are easily expressible with a sufficient degree of mathematical precision: we only need to define precisely what structural properties are preferable; then for a given theory $T$ we can select the signatures $\sigma$ for which the $\sigma$-axiomatization of $T$ has the preferred property. On the other hand it is  less transparent how to give a precise mathematical formulation of the reasons internal to $T$ for choosing a signature. For example, which mathematical criteria allow us to recognize $\bp{\cdot,e}$ as a better signature than $\bp{R,e}$ for the formalization of group theory? For example both give $\Pi_2$-axiomatizations,\footnote{Note that $\forall \vec{x}\,\exists! y\,\psi(\vec{x},y)$ for $\psi$ quantifier free is the $\Pi_2$-sentence\\ 
$\forall \vec{x}\,[\exists y\,\psi(\vec{x},y)\wedge \forall w,z\, (\psi(\vec{x},w)\wedge \psi(\vec{x},z)\rightarrow w=z)]$.} although it is clear that any mathematician would regard the $\bp{\cdot,e}$-axiomatization  as more natural than the one given in $\bp{R,e}$. Can we turn this qualitative preference into a precise mathematical criterion?

\raggedbottom

\section{Simple and complicated concepts for a mathematical theory} \label{sec:simpcomplconc}
 
When dealing with a first order theory $T$ in a signature $\tau$, the $\tau$-formulae naturally correspond to the concepts/notions $T$ intends to formalize. The well-known L\'evy complexity then gives a measure gauged by  $\tau$ of the complexity of a concept/notion of $T$ by analyzing the number and patterns of the quantifiers that appear in a $\tau$-formula which should formalize it. The L\'evy hierarchy of the $\tau$-formulae with respect to $T$-equivalence stratifies the concepts/notions of $T$ as follows:
the basic notions are those formalized by a boolean combination of atomic $\tau$-formulae.  The complexity of a notion then increases according to the number of alternations of $\forall,\exists$ quantifiers that a formula $\phi$ formalizing it displays when expressed in a prenex normal form. In assigning a L\'evy complexity to a notion of $T$, we consider, among all the $\tau$-formulae which are $T$-equivalent and formalize it, those in prenex normal form with the least number of quantifier alternations. In this way we can assign complexity $\Pi_n$, $\Sigma_n$ or $\Delta_n$ to the notions of $T$ that are expressible by means of $\tau$-formulas. More precisely,
$\Pi_0=\Sigma_0=\Delta_0$ is the complexity of  notions formalized by boolean combinations of atomic formulae; $\Pi_{n+1}$ represents the complexity of those notions formalized by a formula of type $\forall \vec{x}\,\psi(\vec{x},\vec{y})$ with $\psi(\vec{x},\vec{y})$ $\Sigma_n$; $\Sigma_{n+1}$ represents the complexity of those notions formalized by a formula of type $\exists\vec{x}\,\phi(\vec{x},\vec{y})$ with $\phi(\vec{x},\vec{y})$ $\Pi_n$, and $\Delta_n$ represents the complexity of the notions whose complexity lies in $\Pi_n\cap\Sigma_n$.
However, Levy complexity for a notion of $T$ does not always match with the intuitive complexity we attribute in practice to it, and is clearly dependent on the signature $\tau$ we use to express it. For example the process of Morleyization (see \cite[Section 3.2]{TENZIE} or Def. \ref{Def:AMCSpec} below) produces a mathematically equivalent axiomatization $T^*$ of a first order $\tau$-theory $T$ in a new signature $\tau^*\supseteq\tau$, such that every $\tau^*$-formula is $T^*$-equivalent to an atomic formula; hence all notions of $T$ become of the same Levy complexity when expressed in signature $\tau^*$ according to the mathematically equivalent theory $T^*$.

To better understand the relevance of the signature in describing the correct logical complexity of the notions of a theory, let us examine the case of the  theory of \emph{commutative semi-rings with no zero-divisors}. 
Consider the signature $\bp{+,\cdot,0,1}$, which is standardly used to axiomatize rings and fields. In this signature we can provide the axioms of commutative semi-rings with no zero-divisors 
by means of the following universal sentences:

\begin{align} \label{ax:commsemiringsno0div}
&\forall x, y\,(x\cdot y=y\cdot x),\\ \nonumber
&\forall x,y,z\, [(x\cdot y)\cdot z=x\cdot(y\cdot z)],\\ \nonumber
&\forall x \,(x\cdot 1=x\wedge 1\cdot x=x),\\ \nonumber
&\\ \nonumber
&\forall x, y\,(x + y=y + x),\\ \nonumber
&\forall x,y,z\, [(x+y)+ z=x+(y+ z)],\\ \nonumber
&\forall y \,(x+ 0=x\wedge 0+ x=x),\\ \nonumber
&\\ \nonumber
&\forall x,y,z\, [(x + y)\cdot z=(x\cdot y)+ (x\cdot z)],\\ \nonumber
&\\ \nonumber
&\forall x, y\,[x\cdot y=0\rightarrow (x=0\vee y=0)].
\end{align}

By further extending the above list we can obtain the theory of \emph{commutative rings with no zero-divisors} (e.g. \emph{integral domains}), 
adding the $\Pi_2$-axiom
\begin{align} \label{ax:commringsno0div}
\forall x\exists y\, (x+y=0),
\end{align}

\noindent and the theory of \emph{fields} 
by further adding the  $\Pi_2$-axiom

\begin{align} \label{ax:fieldsno0div}
\forall x\,[x\neq 0\rightarrow \exists y\, (x\cdot y=1)],
\end{align}

\noindent  and finally the theory of \emph{algebraically closed fields} 
by supplementing the above axioms with the  following $\Pi_2$-sentences\footnote{Notice that any model of Axioms (\ref{ax:commsemiringsno0div}) and (\ref{ax:ACF}) is automatically a field.} for all $n\geq 1$
\begin{align} \label{ax:ACF}
\forall x_0\dots x_n\exists y\, \sum x_i\cdot y^i=0.
\end{align}

A common feature of the above set of $\Pi_2$-axioms is that they assert the existence of solutions to certain basic equations of the theory
expressible in terms of sum, multiplication, $0$, $1$, i.e. the atomic formulae of $\bp{\cdot,+,0,1}$.
Indeed,  any commutative semi-ring $\mathcal{M}$ without zero divisors can be extended to an algebraically closed field simply by extending it with solutions to the basic polynomial equations with parameters in $\mathcal{M}$ (and which may not exist in $\mathcal{M}$). 

In the signature $\bp{+,\cdot,0,1}$ we have as basic operations $+,\cdot$. Using this signature for 
integral domains we can also define the additive inverse, $-$, using the atomic formula $(x+z=0)$, while for fields we can add the multiplicative inverse, ${}^{-1}$, using the boolean combination of atomic formulae $(x=0\wedge z=0)\vee(x\neq 0\wedge x\cdot z=1)$.
Notice that even when subtraction and division are partially defined on a semiring, we can still meaningfully interpret these two operations in it, just by adding two unary function symbols $-,{}^{-1}$ for them and by adopting the convention that the corresponding inverse operations are trivially defined on the non invertible elements; this is captured for example by the axioms 
 \begin{align}\label{ax:multinv}
 \forall x\, [(\exists y\, (x\cdot y=1)\wedge x\cdot x^{-1}=1)\vee(\neg\exists y\, (x\cdot y=1)\wedge x^{-1}=0)]
 \end{align}
 to interpret ${}^{-1}$, and
 \begin{align}\label{ax:plusinv}
 \forall x\, [(\exists y\, (x+ y=0)\wedge x+ (-x)=0)\vee(\neg\exists y\, (x+y=0)\wedge (-x)=0)]
 \end{align}
 to interpret $-$.
 
 It is clear that the class of $\bp{+,\cdot,0,1,-,{}^{-1}}$-structures satisfying axioms 
 (\ref{ax:commsemiringsno0div}), (\ref{ax:multinv}), (\ref{ax:plusinv}) are exactly the commutative semi-rings with no zero-divisors, however, expanding the signatures and introducing the new axioms which correctly interpret the new symbols, we perturbed the notion of morphism. To see this, notice that the inclusion of $\mathbb{N}$ into $\mathbb{Z}$ is a 
 $\bp{+,\cdot,0,1}$-morphism but not a $\bp{+,\cdot,0,1,-,{}^{-1}}$-morphism as $-(2)=0$ when computed in $\mathbb{N}$ seen as a model of (\ref{ax:commsemiringsno0div}), (\ref{ax:multinv}), (\ref{ax:plusinv}), while $-(2)=-2$  when computed in  $\mathbb{Z}$ seen as a model of (\ref{ax:commsemiringsno0div}), (\ref{ax:multinv}), (\ref{ax:plusinv}). Similarly the inclusion of
 $\mathbb{Z}$ into $\mathbb{Q}$ is a $\bp{+,\cdot,0,1,-}$-morphism but not a $\bp{+,\cdot,0,1,-,{}^{-1}}$-morphism.
 
Moreover, we can observe that in the class of integral domains as formalized in signature 
$\bp{+,\cdot,0,1,-}$
subtraction is axiomatized now by the universal sentence $\forall x\,[x+(-x)=0]$ rather than the
$\Pi_2$-sentence (\ref{ax:plusinv}), which defines it in the larger class of commutative semirings.
 Similarly, if we consider the class of 
  $\bp{+,\cdot,0,1,-,{}^{-1}}$-structures which are fields, i.e. models of (\ref{ax:commsemiringsno0div}), 
 (\ref{ax:commringsno0div}), (\ref{ax:fieldsno0div}), then (\ref{ax:multinv}) becomes logically equivalent, modulo the other axioms, to
 $\forall x\,[x\neq 0\rightarrow x\cdot x^{-1}=1]$, which is a universal $\bp{+,\cdot,0,1,-,{}^{-1}}$-sentence.
 On the algebraic side note that the $\bp{\cdot,+,0,1}$-morphisms between integral domains naturally extend to $\bp{\cdot,+,0,1,-}$-morphisms, as the operation $x\mapsto-x$ is preserved by additive morphisms on rings. Similarly the $\bp{\cdot,+,0,1}$-morphisms
  between fields naturally extend to $\bp{+,\cdot,0,1,-,{}^{-1}}$-morphisms, as the operation 
  $x\mapsto x^{-1}$ is preserved by additive and multiplicative morphisms between fields.

More generally, when a $\tau$-theory $T$ can define an operation by means of universal $\tau$-equations, the $\tau$-morphisms between $\tau$-models of $T$ preserve the operation. In our concrete example, $\bp{+,\cdot,0,1}$-morphisms between integral domains naturally extends to $\bp{+,\cdot,0,1,-}$-morphisms, while $\bp{+,\cdot,0,1}$-morphisms between fields naturally extends to $\bp{+,\cdot,0,1,-,^{-1}}$-morphisms exactly because the operations $-,^{-1}$ are defined by universal equations in the respective theories. 

These examples suggest the following observation: when considering the first order axiomatization of a theory, we should carefully consider not only its class of models, but also the class of morphisms between these models.
Dealing with the theory of commutative semi-rings with no zero-divisors, the notion of additive inverse has a complexity which exceeds that of $+$ and $\cdot$. We can detect this by noticing that the class of morphisms between these structures shrinks when we impose the preservation of this operation. On the other hand, in the context of commutative rings, the notion of additive inverse has the same complexity of $+,\cdot$. But, again, this  is not the case for the notion of multiplicative inverse: we need to focus on the theory of field in order to regard it as a concept with the same complexity as addition and multiplication.
We can sum up these considerations as follows:
\begin{itemize}
\item
The signature $\bp{+,\cdot,0,1}$ is suitable for commutative (semi)rings with no zero-divisors, and fields;
\item 
The signature $\bp{+,\cdot,0,1,-}$ is better suited for commutative rings with no zero-divisors, and fields;
\item
The signature $\bp{+,\cdot,0,1,-,{}^{-1}}$ is better suited (only) for fields.
\end{itemize}

The algebraic case we discussed is instructive and guides us towards more general considerations on the interplay between the complexity of concepts/notions of a theory $T$, the first order language $\tau$ in which to express them, and the properties of the $\tau$-morphisms between models of $T$. First of all, the L\'evy hierarchy gives a potentially useful hierarchy by which we can stratify the complexity of concepts/notions of a theory; furthermore our discussion suggests that for a given a signature $\tau$, the basic concepts of a $\tau$-theory $T$ are not only those expressed by (boolean combinations of) atomic $\tau$-formulae, but also those expressed by universal $\tau$-equations (e.g. the additive inverse is a simple operation for ring theory, being axiomatized in signature $\bp{+,\cdot,0,1,-}$ by universal quantification on an equation, while it is not for semi-ring theory as axiomatized in signature $\bp{+,\cdot,0,1,-}$). The reason is that if $\mathcal{M}\sqsubseteq\mathcal{N}$ are $\tau$-structures which model the same $\tau$-theory $T$, the graph of a relation/operation $R^{\mathcal{M}}$ defined on $\mathcal{M}$ as prescribed by the universal axiom $\psi_R$ of $T$ is exactly the restriction to $\mathcal{M}$ of the graph of that same relation/operation $R^{\mathcal{N}}$ defined now on $\mathcal{N}$ as prescribed by $\psi_R$ (in our example, $-$ is preserved by $\bp{+,\cdot,0,1,-}$-morphisms of rings but not by $\bp{+,\cdot,0,1,-}$-morphisms of semi-rings).

Summing up, when dealing with a mathematical theory $T$ for which we have a clear picture of its intended class of models, the choice of a signature $\tau$ in which to formalize $T$ can be made on the basis of which morphism between models of $T$ we would like to be the $\tau$-morphisms between the $\tau$-models of $T$, or, correspondingly, on the basis of which concepts of $T$ we would like to be of low logical complexity in the L\'evy hierarchy induced by $\tau$ on $T$.  Once this decision is made, this affects the type of signature $\tau$ in which the theory can be formalized, what are the $\tau$-morphisms between the $T$-models, and what are the universal $\tau$-axioms of $T$.
Therefore, we can get different stratifications of the complexity of concepts for a theory $T$ by selecting the appropriate class of morphisms between its models and  by selecting which concepts of $T$ we require to be axiomatized by universal sentences.

\section{Existentially closed models, model completeness, and model companionship} \label{sec:exclosmod}

It is now the time to introduce the notion of existentially closed model and of model companionship.\footnote{Reference texts for existentially closed models, model completeness, and model companionship are \cite{CHAKEI90} and the second author's notes \cite{VIAMODCOMPNOTES}.} These model-theoretic notions will be the structural properties of first order theories we will use to select
good signatures for set theory.
Before giving the precise mathematical definitions of these notions, let us start by introducing these concepts in the context of (semi-)ring theory.

A common trait of the set of axioms 
(\ref{ax:commringsno0div}), (\ref{ax:fieldsno0div}), (\ref{ax:ACF}) is that they assert the existence of solutions to certain basic equations that are expressible by atomic formulae in the $\bp{+,\cdot,0,1}$-theory given by axioms (\ref{ax:commsemiringsno0div}).
Indeed,  any commutative (semi-)ring $\mathcal{M}$ without zero divisors can be enlarged to an algebraically closed field simply by extending it with solutions to the basic polynomial equations with parameters in $\mathcal{M}$ (which in principle may not exist in $\mathcal{M}$). 
Notice that the process of adding new solutions to a commutative (semi-)ring can be seen as closing it with respect to its basic operations. It is here  that the choice of a signature $\tau$  for the formalization of a  theory $T$ becomes relevant: the richer $\tau$ is, the more complex are the concepts expressible by atomic $\tau$-formulae (expressing the basic operations and relations) and universal $\tau$-sentences (expressing the basic properties of the $T$-definable relations and operations). Therefore, the richer $\tau$ is, the more closed-off are the models of a universal $\tau$-theory $T$ that have solutions for the atomic $\tau$-formulae with parameters in the model.

This closing-off process is one of the driving forces of mathematics. For example it is by adding new solutions to basic equations that commutative semi-rings with no zero-divisors (like $\mathbb{N}$) have been extended to  commutative rings (like $\mathbb{Z}$, which contains all additive inverses), and then to fields (like $\mathbb{Q}$, which contains all  multiplicative inverses), and finally to algebraically closed fields (like $\mathbb{C}$, in which all diophantine equations have a solution).

Of course, the choice of the signature $\tau$ in which a theory $T$ is formalized affects the outcome of this closing-off process. A basic request is that the $\tau$-theory $T'$ describing the validities in the closed-off structures, should agree with $T$ at least with respect to the $\Pi_1$-consequences of $T$ in signature $\tau$ (as the latter express the properties of the basic concepts and operations of $T$ according to $\tau$). Indeed, only in this case we can say that the closing-off process has been performed with respect to the basic concepts and operations of $T$. 
The notion of existentially closed model defines precisely which structures are the outcome of this closing-off process.

\begin{definition}\label{ecM}
Let $\tau$ be a first order signature and $T$ be a $\tau$-theory.
A $\tau$-structure $\mathcal{M}$ is $T$-ec if:
\begin{itemize}
\item There is some $\mathcal{N}$ $\tau$-superstructure of $\mathcal{M}$ which models $T$.
\item $\mathcal{M}$ is a $\Sigma_1$-substructure of $\mathcal{P}$ whenever $\mathcal{P}$ is a $\tau$-superstructure of $\mathcal{M}$ which models $T$.
\end{itemize}
\end{definition}

A caveat is in order: while a mathematical theory $T$ can be unambiguosly defined regardless of any of its possible first order axiomatizations, to define unambiguously the notion of $T$-ec model we commit ourselves to choosing a specific signature $\tau$ in which to formalize $T$.

\begin{notation}
Given a signature $\tau$ and a $\tau$-theory $T$, $T_\forall^\tau$ is the set of universal $\tau$-sentences which follows from $T$ (similarly we define $T^\tau_{\exists}$,  $T^\tau_{\forall\exists}$, etc.).
\end{notation}

While the set of logical consequences of $T$ is independent of the signature in which we axiomatize it, 
the set $T_\forall^\tau$ heavily depends  on the choice of $\tau$. 

Notice that any $T$-ec $\tau$-structure is a model of $T^\tau_\forall$.
A less trivial observation is the following.

\begin{fact}\cite[Fact 8]{VIAMODCOMPNOTES}
Given a $\tau$-theory $T$, a $\tau$-structure $\mathcal{M}$ is $T$-ec if and only if it is 
 $T^\tau_\forall$-ec. Hence a $\tau$-structure $\mathcal{M}$ which is $T$-ec is also $S$-ec for any $\tau$-theory $S$ such that $S^\tau_\forall=T^\tau_\forall$.
 \end{fact}

It is not hard to see that algebraically closed fields, in the signature $\sigma=\bp{+,\cdot,0,1}$, are $S$-ec structures for $S$ the $\sigma$-theory of commutative
(semi-)rings with no zero-divisors.
Consider the two $\sigma$-theories $U$ and $R$, respectively, of fields and of integral domains that are not fields. Then,  any $\sigma$-structure  $\mathcal{M}$ which is an algebraically closed field is automatically $S$-ec, $U$-ec, and $R$-ec at the same time (since $R^\sigma_\forall=S^\sigma_\forall=U^\sigma_\forall$).

In general, given a $\tau$-theory $T$, the class of $T$-ec models might not be elementary. In the specific case of algebraically closed fields this is the case, since this class of structures is axiomatized by the $\bp{+,\cdot,0,1}$-theory $\bool{ACF}$ given by axioms (\ref{ax:commsemiringsno0div}) and (\ref{ax:ACF}); but this is a rather peculiar case, which also depends on the  specific choice of the signature. 
On the contrary, if we consider the  $\bp{\cdot,e,{}^{-1}}$-theory $T$ of groups, then the $T$-ec $\bp{\cdot,e,{}^{-1}}$-structures do not form an elementary class\footnote{Moreover, none of the signatures for group theory that we presented in Section \ref{sec:rightsignmathth} have existentially closed models which give rise to an elementary class.}  \cite[Example 3.5.16]{CHAKEI90}. Abraham Robinson introduced the notion of model companionship to describe exactly when the collection of $T$-ec $\tau$-models form an elementary class:

 \begin{definition}\label{Def:MC}
Let $\tau$ be a signature.
\begin{itemize}
\item A $\tau$-theory $T$ is \emph{model complete}\footnote{Usually model completeness is defined by requiring that the substructure relation coincides with the elementary substructure relation \cite[Def. on page 186]{CHAKEI90}. Our definition is equivalent in view of \cite[Prop. 3.5.15]{CHAKEI90} combined with \cite[Lemma 2.1.1]{viale2021absolute}.} if for any
$\tau$-structure $\mathcal{M}$ 
\[
\mathcal{M}\models T \text{ if and only if $\mathcal{M}$ is $T^\tau_\forall$-ec}.
\]
\item 
A $\tau$-theory $T$ is the \emph{model companion} of a $\tau$-theory $S$ if:
\begin{enumerate}
\item[(i)] $S$ and $T$ are jointly consistent: i.e. $T^\tau_\forall=S^\tau_\forall$,
(or equivalently ---by \cite[Lemma 2.1.1]{viale2021absolute} or  \cite[Remark 3.5.6(2)]{CHAKEI90}--- every model of $T$ embeds in a model of $S$, and vice versa);
\item[(ii)] $T$ is model complete. 
\end{enumerate}
\end{itemize}
\end{definition}

It is useful to recall that model complete theories are axiomatized by their 
$\Pi_2$-consequences (see \cite[Prop. 3.5.10]{CHAKEI90}), and that the model companion of a $\tau$-theory $T$ (if it exists) is unique (by \cite[Prop. 3.5.13]{CHAKEI90}).

%
%

Conforming to our previous notation of $R$ being the $\bp{+,\cdot,0,1}$-theory of integral domains which are not fields, and $U$ being the $\bp{+,\cdot,0,1}$-theory of fields, we can observe that no algebraically closed field is a model of $R$, and, conversely, no model of $R$ is an algebraically closed field. As a matter of fact, the notion of model companionship for a $\tau$-theory $T$ does not necessarily isolate a proper subclass of the models of $T$, but rather the models that are closed-off with respect to the basic operations and relations definable in $T$ by universal $\tau$-sentences.

%
%

%


Given the connection between model companions and existentially closed structures, it should be no surprise that the notion of model companionship is dependent on the choice of the signature. For example: consider the language 
$\bp{+,\cdot,0,1,{}^{-1}}$ and  the theory $S'$ which results from adding axiom (\ref{ax:multinv}) to axioms (\ref{ax:commsemiringsno0div}), (\ref{ax:commringsno0div}).
 We obtain that $S'$ is still the theory of commutative rings with no zero-divisors, now formalized in the signature  $\bp{+,\cdot,0,1,{}^{-1}}$. 
  However,  the  $\bp{+,\cdot,0,1,{}^{-1}}$-theory of algebraically closed fields $\bool{ACF}'$ given by axioms (\ref{ax:commsemiringsno0div}), (\ref{ax:ACF}), (\ref{ax:multinv}) is not  the model companion of the  $\bp{+,\cdot,0,1,{}^{-1}}$-theory $S'$. In order to see this, let
 $\mathcal{M}$ be the field of complex numbers and let $\mathcal{N}$ be $\mathbb{C}[X]$ (the ring of polynomials in complex coefficients and variable $X$), both seen as  $\bp{+,\cdot,0,1,{}^{-1}}$-structures. Then it is the case that $\mathcal{M}\sqsubseteq\mathcal{N}$, but now $\mathcal{N}$ models the $\Sigma_1$-sentence for $\bp{+,\cdot,0,1,{}^{-1}}$
 $\exists y\,(y\neq 0\wedge y^{-1}=0)$ (as witnessed by the polynomial $X$), while $\mathcal{M}$ does not.
 In particular, $\mathcal{M}$ is not $S'$-ec, as $\mathcal{M}$ is not a $\Sigma_1$-substructure of 
 $\mathcal{N}$, which is a model of  $S'$.
 
Given that the existence of a(n absolute) model companion for a theory $T$ depends so much on the choice of the signature, we can use the former to select the latter: 

\medskip

\begin{description}
\item[Model-theoretic criterion for selecting signatures] \label{modthcritforsignST}
Given a mathematical theory $T$ and two signatures $\tau$ and $\sigma$ in which $T$ can be axiomatized, $\tau$ is preferable to $\sigma$ for $T$ if the $\tau$-axiomatization of $T$ admits a(n absolute) model companion, while its $\sigma$-axiomatization does not.
\end{description}

\medskip
Clearly this criterion is in the mold of \ref{ABSJUST}.
Let us see with a concrete example why criteria of this kind can cohere with others that are internal to a theory.
Recall the example of integral domains: in Section \ref{sec:simpcomplconc} we argued on the basis of considerations internal to this theory that the signature $\bp{\cdot,+,0,1,-}$ is better suited for it than $\bp{\cdot,+,0,1,-,^{-1}}$.
Let us apply the above abstract criteria to select which among the two signatures is best suited for the theory of integral domains. Once again the answer is $\bp{\cdot,+,0,1,-}$, in this case because in this signature the theory of integral domains has a natural \emph{and simple to compute} model companion (which is the $\bp{\cdot,+,0,1,-}$-theory of algebraically closed fields).\footnote{On the other hand (given our limited knowledge of algebra) a model companionship result for the theory of integral domains in signature $\bp{\cdot,+,0,1,-,^{-1}}$ is a much harder problem (of which we don't know the answer); what we can say for sure is that the eventual model companion of integral domains for $\bp{\cdot,+,0,1,-,^{-1}}$ won't be the theory of algebraically closed fields, and that the theory of integral domains does not have an absolute model companion in the above signature.}
It is actually the case (at least for concrete examples) that  we can prove that some theory $T$ have a(n absolute) model companionin in a signature $\tau$ only when we outline reasons internal to $T$ of why the universal $\tau$-theory of $T$ is catching certain central and basic concepts of $T$.

A category theoretic perspective can also be useful to justify this model-theoretic criterion. We already remarked that for a mathematical theory $T$, its semantic extension might be given not only by its class of models, but also by its intended class of  \emph{morphisms between these models}, i.e. by a certain category $\mathcal{E}_T$ of models and morphisms (this is the standard approach with which ring theory is presented in a basic undergraduate course in algebra).
By axiomatizing $T$ in some signature $\tau$ we univocally describe the objects of $\mathcal{E}_T$ as the $\tau$-models of $T$, but the arrows of $\mathcal{E}_T$ may not coincide with the $\tau$-morphisms between $\tau$-models of $T$ (consider again the case of integral domains with the two signatures $\bp{\cdot,+,0,1,-}$, $\bp{\cdot,+,0,1,-,^{-1}}$). 
How do we detect the most natural family of morphisms between models of a mathematical theory $T$ (i.e. the arrows of $\mathcal{E}_T$)? Can we even say that there is always a most natural family of morphisms between models of $T$ (e.g. that $\mathcal{E}_T$ can be meaningfully and univocally defined for any mathematical theory $T$)?
Our model theoretic criterion to select signatures is a useful tool to address these questions.

\section{Partial Morleyizations, absolute model companionship, and the AMC-spectrum of a theory} \label{sec:partmorlAMCspec}

We can now give a more formal treatment of the ideas presented in the previous sections, by introducing the central concepts of (partial) Morleyization and 
absolute model companionship.

\subsection{Partial Morleyizations}

\begin{notation}\label{not:keynotation0}
Given a signature $\tau$, 
let $\phi(x_0,\dots,x_n)$ be a $\tau$-formula.
  
We let:
\begin{itemize} 
\item
$R_\phi$ be
a new
$n+1$-ary relation symbols,
\item
$f_\phi$ be a
new $n$-ary function symbols\footnote{As usual we won't distinguish between $0$-ary function symbols with constants.} 
\item
$c_\tau$ be a new
constant symbol.
\end{itemize}
 We also let:
\[
\AX^0_\phi:= \forall\vec{x}[\phi(\vec{x})\leftrightarrow R_\phi(\vec{x})],
\]
\begin{align*}
\AX^1_\phi:= &\forall x_1,\dots,x_n\,\\
&[(\exists!  y\phi(y,x_1,\dots,x_n)\rightarrow \phi(f_\phi(x_1,\dots,x_n),x_1,\dots,x_n))\wedge\\
&\wedge(\neg\exists!  y\phi(y,x_1,\dots,x_n)\rightarrow f_\phi(x_1,\dots,x_n)=c_\tau)]
\end{align*}
for $\phi(x_0,\dots,x_n)$ having at least two free variables, and
\[
\AX^1_\phi:= \qp{(\exists!  y\phi(y))\rightarrow \phi(f_\phi)}\wedge
 \qp{(\neg\exists!  y\phi(y))\rightarrow c_\tau=f_\phi}.
\]
for $\phi(x)$ having exactly one free variable.

Let $\bool{Form}_\tau$ denotes the set of $\tau$-formulae.
For $A\subseteq \bool{Form}_\tau\times 2$
 \begin{itemize}
 \item 
 $\tau_A$ is the signature obtained by adding to $\tau$ all relation symbols $R_\phi$, for
 $(\phi,0)\in A$, and all function symbols  $f_\phi$, for $(\phi,1)\in A$ (together with the special symbol $c_\tau$ if at least one $(\phi,1)$ is in $A$).
 \item
  $T_{\tau,A}$ is the $\tau_A$-theory having as axioms
the sentences $\AX^i_\phi$ for $(\phi,i)\in A$.
\end{itemize}
We let $\tau^*=\tau_C$ for $C=\bool{Form}_\tau\times \bp{0}$ and $T^*$ be the $\tau^*$-theory $T_{\tau,C}$.

\end{notation}

Roughly speaking the axioms $\AX^i_\phi$ are meant to formalize the idea that the (potentially complicated) concept formalized by $\phi$ according to $\tau$ becomes an elementary concept if we formalize instead it using $\tau_{\bp{(\phi,i)}}$. More precisely, given a $\tau$-theory $S$:
\begin{itemize}
\item
$\AX^0_\phi$ transforms the concept of $S$ formalized by the $\tau$-formula $\phi$ (which could have a high Levy complexity according to $\tau$) in  a concept formalized by an atomic $\tau\cup\bp{R_\phi}$-formula according to $S+\AX^0_\phi$.  
\item For any $\tau$-theory $S$, and for $C=\bool{Form}_\tau\times \bp{0}$,
$S^*=S+T_{\tau,C}$ is a $\tau^*=\tau_C$-theory admitting quantifier elimination (the Morleyization of $T$).\footnote{See \cite[Section 3.2, pp. 31-32]{TENZIE}.}
\item
$\AX^1_\phi$ performs a delicate Skolemization process. The standard Skolemization introduces for each formula $\phi$ a Skolem term $f_\phi$ together with an axiom stating
$\forall\vec{x}\,[\exists y\,\phi(\vec{x},y)\leftrightarrow \phi(\vec{x},f_\phi(\vec{x}))]$. 
However this has serious drawbacks. Consider the following situation: we have a $\tau$-structure $\mathcal{M}$ which is such that for all $b\in \mathcal{M}$ $\exists y\,\phi(b,y)$ holds, and also such that for some $a\in \mathcal{M}$, we have that
$\phi(a,y)$ has more than one solution in $\mathcal{M}$. In this case we can find distinct functions 
$h_0,h_1:\mathcal{M}\to\mathcal{M}$ such that $\forall x\,[\exists y\,\phi(x,y)\leftrightarrow \phi(x,f_\phi(x))]$ is satisfied in the expansions of $\mathcal{M}$ to a $\tau\cup\bp{f_\phi}$-structure which interprets $f_\phi$ by either one of the $h_i$. Which one of the two expansions should we consider the natural Skolemization of $\mathcal{M}$ for the formula $\phi$? $\AX^1_\phi$ rules out this source of ambiguity by saying that when $\phi(a,y)$ has more than one solution in $\mathcal{M}$, then $f_\phi^\mathcal{M}$ does not select any of them and instead behaves trivially by assigning to $a$ the value $c_\tau^\mathcal{M}$. In this way, given a formula $\psi(\vec{x},y)$, once we choose the interpretation of $c_\tau$, any $\tau$-structure has a unique expansion to a $\tau\cup\bp{f_\psi,c_\tau}$-structure which is a model of $\AX^1_\psi$ and has that interpretation of $c_\tau$.
\end{itemize}

\smallskip

Notice that with respect to the signatures expanding $\bp{\cdot,+,0,1}$ discussed in the previous section, Axiom (\ref{ax:multinv}) is  (equivalent to) $\AX^1_\phi$ for
$\phi(x,y)$ being $x\cdot y=1$ (assuming $c_\tau$ is interpreted by $0$). Similarly Axiom
(\ref{ax:plusinv}) is (equivalent to) $\AX^1_\psi$ for
$\psi(x,y)$ being $x+ y=0$.
In what follows we want to analyze what happens when the Morleyization process is performed on arbitrary subsets of $\bool{Form}_\tau\times 2$.

\subsection{Absolute model companionship}\label{subsec:AMC}

The following properties will bring us to introduce the notion of absolute model companionship. 

\begin{itemize}
\item
A $\tau$-structure $\mathcal{M}$ is $T$-ec if and only if it is $T^\tau_\forall$-ec (by \cite[Fact 8]{VIAMODCOMPNOTES}).
\item
If $T$ is \emph{complete},
a $T$-ec structure $\mathcal{M}$ realizes any $\Pi_2$-sentence
which is \emph{consistent} with $T^\tau_\forall$ (by \cite[Fact 12]{VIAMODCOMPNOTES}).
\item
If a $\tau$-theory $S$ is the model companion of a $\tau$-theory $T$, $S$ is axiomatized by its $\Pi_2$-consequences for signature $\tau$ (by Robinson's test, e.g. \cite[Lemma 14]{VIAMODCOMPNOTES}).
\end{itemize}

These results are outlining ``syntactic properties'' of $T$-ec models which reflect the ``semantic fact'' that $T$-ec models are the ``closed-off'' structures for $T$; furthermore they outline that the axiomatization of the model companion of a \emph{complete} theory can be detected by separately solving a consistency task for each $\Pi_2$-sentence for $\tau$.  One might therefore wonder whether these syntactic properties can give another equivalent characterization of model companionship, i.e. if the model companion $S$ of a $\tau$-theory $T$ (whenever it exists) can be axiomatized by the family of all $\Pi_2$-sentences for $\tau$ which holds in some model of the universal and existential fragment of $T$. This is not always the case; it holds when $T$ is a \emph{complete} model companionable theory, however there is a standard counterexample: the $\bp{+,\cdot,0,1}$-theory $\bool{Fields}_0$ of fields of characteristic $0$, given by axioms (\ref{ax:commsemiringsno0div}), (\ref{ax:commringsno0div}), (\ref{ax:fieldsno0div}) and the infinitely many axioms granting that the characteristic of its models is $0$ (note that $\bool{Fields}_0$ is not a complete theory) has as its model companion $\bool{ACF}_0$, which adds axioms (\ref{ax:ACF})  to $\bool{Fields}_0$ (note that $\bool{ACF}_0$ is a complete theory). Indeed, the sentence $\forall x \neg(x^2+1=0)$ is a $\Pi_2$-sentence (actually $\Pi_1$) which is not in $\bool{ACF}_0$ but holds in $\mathbb{Q}$, hence is consistent with the universal and existential fragment of $\bool{Fields}_0$. 

Absolute model companionship rules out counterexamples of this kind to the above tentative characterization of model companionship. 

\begin{notation}
Given a $\tau$-theory $T$,
$T_{\forall\vee\exists}^\tau$ denotes the \emph{boolean combinations} of universal $\tau$-sentences which follow from $\tau$. 
\end{notation}
Note that $T_{\forall\vee\exists}^\tau$ may convey more information than $T_{\forall}^\tau\cup T_{\exists}^\tau$, as there could be $\psi\vee\phi$ (with $\psi$ universal and $\phi$ existential) which is in the former but does not follow from the latter.

\begin{definition}\cite[Def. 2.2.4, Lemma 2.3.4]{viale2021absolute} \label{def:AMC}
Given a $\tau$-theory $T$, we say that a $\tau$-sentence $\psi$ is \emph{strongly $T_{\forall\vee\exists}^\tau$-consistent for $\tau$} whenever
$\psi+R^\tau_{\forall\vee\exists}$ is consistent for any $\tau$-theory $R$ which extends $T$.

A $\tau$-theory $S$ is the \emph{absolute model companion} (AMC) of a $\tau$-theory $T$ if it is the model companion of $T$ and it is axiomatized by the $\Pi_2$-sentences (for $\tau$) which are strongly $T_{\forall\vee\exists}^\tau$-consistent for $\tau$.

Equivalently (by \cite[Lemma 22]{VIAMODCOMPNOTES}) $T$ is the AMC of $S$ if:
\begin{itemize}
\item $T^\tau_{\forall\vee\exists}=S^\tau_{\forall\vee\exists}$,
\item $T$ is model complete.
\end{itemize}
\end{definition}

Note that a model complete $\tau$-theory $T$ is the model companion of $T^\tau_\forall$ and the AMC of 
$T^\tau_{\forall\vee\exists}$.

A $\tau$-theory $S$ can have an AMC if, to some extent, it has already maximized the family of existential sentences which are consistent with its universal consequences. 
Let us elaborate more on this point.
Given a $\tau$-theory $T$, the family of strongly $T_{\forall\vee\exists}^\tau$-consistent $\Pi_2$-sentences for $\tau$ describes a fragment of the $\Pi_2$-sentences which hold in a $T$-ec $\tau$-structure (see \cite[Fact 2.2.5]{viale2021absolute}). An AMC describes the set-up in which the choice of $\tau$ is such that the family of strongly $T_{\forall\vee\exists}^\tau$-consistent $\Pi_2$-sentences  axiomatizes a model complete theory.
This leads to another peculiar aspect separating model-companionship from absolute model companionship:
\begin{lemma}\cite[Lemma 1.5]{viale2021absolute} \label{lem:charAMCvsMC}
Assume $S,S'$ are $\tau$-theories such that $S'$ is the AMC of $S$.
Then for all $T\supseteq S$, $S'+T_{\forall}$ is the AMC of $T$.
\end{lemma}
On the contrary this property fails for model companionship: the universal fragment of the theory of $(\mathbb{Q},+,\cdot,0,1)$ is inconsistent with the $\bp{+,\cdot,0,1}$-theory of algebraically closed fields.
As we already mentioned, complete first order $\tau$-theories are
model companionable if and only if they admit an AMC. Indeed, for a complete theory $S$, the information conveyed by $S^\tau_\forall$ and by $S^\tau_{\forall\vee\exists}$ is the same. On the other hand we have already seen that there are non-complete $\tau$-theories admitting a model companion but not an AMC, e.g. the $\bp{+,\cdot,0,1}$-theory of fields.
Let us briefly analyze more in details what is the obstruction for the $\bp{+,\cdot,0,1}$-theory of fields of characteristic $0$ to have the $\bp{+,\cdot,0,1}$-theory of algebraically closed fields of characteristic $0$ as its AMC. The problem is that 
Diophantine equations can be expressed by atomic formulae of $\bp{+,\cdot,0,1}$ and not all of  them have rational solutions. On the other hand the existence of solutions for these equations does not contradict the field axioms. Note that the smallest existentially closed $\bp{+,\cdot,0,1}$-model 
$\mathcal{M}$ for the theory of fields is given by the algebraically closed field consisting exactly of the solutions of diophantine equations. Its standard construction builds it as a direct limit of Galois extensions of the rationals each adding new solutions to polynomial equations with coefficients in $\mathbb{N}$ but without rational solutions. By performing this construction we can preserve the field axioms and the characteristic $0$ (which holds in $\mathbb{Q}$ as well as in all the Galois extensions of $\mathbb{Q}$ under consideration), but at each stage we may invalidate some universal sentence that is true in $\mathbb{Q}$ (i.e. a universal sentence asserting that a certain diophantine equation does not have rational solutions). In particular the closing-off process performed in signature $\sigma=\bp{+,\cdot,0,1}$ which brings from the theory of fields $T$ to that of algebraically closed fields cannot be described only by the $\Pi_2$-sentences $\psi$ which are strongly $T_{\forall\vee\exists}^\tau$-consistent for $\sigma$. Indeed, there are other $\Pi_2$-sentences which need to be realized in an algebraically closed field and which are not strongly $T_{\forall\vee\exists}^\tau$-consistent for $\sigma$
(i.e. the existential statements asserting the existence of solutions for diophantine equations). 
In particular the $\bp{+,\cdot,0,1}$-theory of algebraically closed fields is the AMC only of the theory of integral domains which extends the field of algebraic numbers (i.e. which have solutions to all polynomial equations expressible in the signature $\bp{+,\cdot,0,1}$).

On the other hand we will argue that in the case of set theory, AMC faithfully describes the closing-off process of models of set theory with respect to basic set-theoretic operations.

\subsection{The (A)MC-spectra of a first order theory}

Model theory has been extremely successful in classifying the complexity of a mathematical theory according to its structural properties, and has produced a variety of properties and criteria to separate the (so called) tame mathematical theories from the (so called) untame ones (e.g. $o$-minimality, stability, simplicity, NIP), see \cite{baldwin} for a philosophical appraisal of this theme. Typically a mathematical theory is considered hard to classify (and thus untame) if it can code in itself first order arithmetic. In this respect the $\in$-theory $\ZFC$ is untame.

We have already observed that most mathematical theories admit many different first order axiomatizations in (almost) as many distinct signatures.  A common characteristic of tameness properties such as $o$-minimality, stability, simplicity, and NIP is that they are \emph{signature invariant}. More precisely, if we take a $\tau$-theory $T$ and we consider its Morleyization $T^*$ in the signature $\tau^*$ (see Notation \ref{not:keynotation0}), $T$ is $o$-minimal (stable, simple, NIP) if and only if so is $T^*$. In contrast, this is not the case for Robinson's notion of model companionship, as we have already observed for the theory of integral domains.
We instantiate the model-theoretic criterion to select signatures formulated in  the last part of Section  \ref{sec:exclosmod} by means of the following: 

\begin{definition} \label{Def:AMCSpec}
Let $T$ be a $\tau$-theory.
\begin{itemize}
\item
Its  AMC-spectrum ($\SpecAMC{T}{\tau}$) is given by the sets $A\subseteq \bool{Form}_\tau\times 2$ such that $T+T_{\tau,A}$ has an AMC (which we denote by $\bool{AMC}(T,A)$).
\item
Its MC-spectrum  ($\mathfrak{spec}_{\bool{MC}}(T,\tau)$)  is given by the sets $A\subseteq \bool{Form}_\tau\times 2$ such that $T+T_{\tau,A}$ has a model companion  (which we denote by $\bool{MC}(T,A)$).
\end{itemize}
\end{definition}

Clearly $\mathfrak{spec}_{\bool{MC}}(T,\tau)$ is a superset of $\SpecAMC{T}{\tau}$ and the two can be distinct.
Observe that $C=\bool{Form}_\tau\times \bp{0}$ is always in the AMC-spectrum of a theory $T$, as $T+T_{\tau,C}$ admits quantifier elimination and therefore it is model complete and its own AMC in signature $\tau_C=\tau^*$. Moreover, $\emptyset $ is in the (A)MC-spectrum of $T$ if and only if $T$ has an (absolute) model companion.
For $\sigma=\bp{+,\cdot,0,1}$ and the various $\sigma$-theories considered in Section \ref{sec:simpcomplconc} we have that:
\begin{itemize} 
\item
$\emptyset$ is in 
$\mathfrak{spec}_{\bool{MC}}(T,\sigma)\setminus\mathfrak{spec}_{\bool{AMC}}(T,\sigma)$ for $T$ the
$\sigma$-theory of commutative (semi)rings with no $0$-divisors or the theory of fields.
\item 
For $\phi$ the $\sigma$-formula defining the additive inverse, 
$\bp{(\phi,1)}$ is in $\mathfrak{spec}_{\bool{MC}}(T,\sigma)$ if $T$ is the theory of commutative rings with no $0$-divisors or the theory of fields, while
$\bp{(\phi,1)}$ is not in 
$\mathfrak{spec}_{\bool{MC}}(S,\sigma)$ if $S$ is the $\sigma$-theory of commutative semirings with no $0$-divisors. 
\item  
For $\psi$ the $\sigma$-formula defining the multiplicative inverse, 
$\bp{(\psi,1)}$ is in $\mathfrak{spec}_{\bool{MC}}(T,\sigma)$ if $T$ is the theory of fields, while
$\bp{(\psi,1)}$ is not in 
$\mathfrak{spec}_{\bool{MC}}(S,\sigma)$ if $S$ is the $\sigma$-theory of commutative (semi)rings with no $0$-divisors. 
\item
It can be shown that $\mathfrak{spec}_{\bool{AMC}}(T,\sigma)=\emptyset$ for  $T$ the $\sigma$-theory of commutative (semi)rings with no $0$-divisors or the theory of fields.
\end{itemize}


A category theoretic perspective on the notion of (A)MC spectrum can again be enlightening: given a $\tau$-theory $T$, consider the category $\mathcal{C}_T$ given by the elementary class of its $\tau$-models as objects and the $\tau$-morphisms (e.g. maps which preserve the \emph{atomic} $\tau$-formulae) between them as arrows. By taking a subset $A$ of $\bool{Form}_\tau\times 2$ we can pass to the category $\mathcal{C}_{T+T_{\tau,A}}$ whose objects are
 $\tau_A$-models of $T+T_{\tau,A}$ and whose arrows are the $\tau_A$-morphisms (in accordance with Notation \ref{not:keynotation0}). There is a natural identification of the objects of $\mathcal{C}_T$ and those of $\mathcal{C}_{T+T_{\tau,A}}$, but the arrows of $\mathcal{C}_{T+T_{\tau,A}}$ are now a possibly much narrower subfamily of the arrows of $\mathcal{C}_{T}$. Important structural information on a $\tau$-theory $T$ is given (at least from a categorial point of view) by classifying which sets $A$ produce an elementary class of 
 $\tau_A$-models for which the $T+T_{\tau,A}$-ec models constitute themselves an elementary class.
 The model companionship spectrum of $T$ gives exactly this structural information on the arrows of $\mathcal{C}_T$.
 Similar considerations apply for the AMC-spectrum of a theory. 
 On the other hand simplicity, stability, NIP
 provide fundamental structural informations on the class of objects of $\mathcal{C}_T$, but it is not transparent whether they also convey information on its class of arrows.

\section{The (A)MC-spectra of set theory} \label{sec:AMCspecST}

From now on we focus our attention on set theory. Our main goal is to study the model companions of set theory and the relevant signatures for which they exist. We will not only provide a new set of arguments in favour of $2^{\aleph_0}=\aleph_2$, but will also present a new method to detect set theoretic validities. The realist-minded readers can consider the following part of the paper as an attempt to produce a complete (first order) axiomatization of set theory. 
Although a realist stance towards mathematics is not needed to understand and prove the results we are going to present, it is within a realist agenda that we can better appreciate their meaning. 
We can make this perspective more explicit by outlining a principle and a goal that will motivate the present model-theoretic approach.

\begin{description}
    \item[Semantic realism] The universe of all sets $(V,\in)$ is sufficiently well-defined  so that its $\in$-theory is a well-defined (complete and non-recursive) first order theory. Therefore, any set theoretical statement axiomatizable in first order logic receives, in $(V,\in)$, a well defined truth value.
    \item[Hilbertian completeness] $(V,\in)$ should offer a \emph{complete} picture of what exists in set theory. That is to say, $(V,\in)$ should contain any set whose existence is not in direct contradiction with already recognized truths of set theory (e.g. $\ZFC$).\footnote{Taken at face value, this seems like an impossible task. As a matter of fact, we can consistently add to $\ZFC$ the statement that there is a bijection between $2^{\aleph_0}$ and $\aleph_1$, but also the statement that the continuum has size $\aleph_2$. Clearly, the two statements are jointly inconsistent. For this reason, we interpret Hilbertian completeness as a motivating goal that needs further refinement in order to be consistently implemented.}
\end{description}


In the context of a realist conception of sets, a consistent formalization of 
Hilbertian completeness would come along with a notion of \emph{maximality} for set theory. Notice that the quest for maximality is what motivates many set-theoretical investigations; for this reason we believe it is appropriate to  analyze the role that model companionship can play in its investigation; however we defer this task to the  the last section of this paper.
Meanwhile, we can propose to implement Hilbertian completeness by means of a more concrete (but still partially informal) principle clearly inspired by it. 

\smallskip

\begin{enumerate}[(IEC)]
\item \label{Descr:IWF}
\textbf{Informal Existential Completeness}
If $P(x,Y)$ is a ``simple'' property formalizable in the $\in$-signature in parameter some set $Y$ of $V$ and $\exists x\,P(x,Y)$ is consistent with the basic principles of set theory, 
then in $V$ a witness of this existential statement should exist.
\end{enumerate}

\smallskip

Our task now consists in spelling out a precise mathematical definition of which set theoretic properties are to be regarded as ``simple'' and which sets $Y$ of $V$ we can accept as parameters for the property $P$ (also attention must be paid to singling out what we exactly mean by ``consistent with the basic principles of set theory'' to avoid trivial counterexamples to \ref{Descr:IWF}).
Towards this end we can take advantage of our previous discussion about signatures and AMC.

\begin{enumerate}[(a)]
\item\label{cond:simplstprop} 
We determine that a property of sets $P$ is ``simple'' by analyzing the peculiarity of set theory as a mathematical theory. After choosing $A\subseteq\bool{Form}_{\bp{\in}}\times 2$ and performing  a partial Morleyization with respect to $A$, we can check whether $P(x,y)$ can be formalized by a quantifier free formula in signature $\bp{\in}_A$. If this is the case we consider $P(x,y)$ simple (with respect to $A$). 
We clearly need criteria for the choice of $A$. Indeed, if $A=\bool{Form}_{\bp{\in}}\times 2$, then any $\in$-formalizable set theoretic property becomes simple, since it is expressible by an atomic formula of $\bp{\in}_A$. In order to avoid trivialities, we use set theoretic practice and forcing as a guiding tool to select the appropriate sets $A$ and, consequently, the simple set theoretic properties. We regard a property $P(\vec{x})$ as simple if one or more of the following conditions are met:
\begin{enumerate}[(i)]
\item \label{cond:simplstprop(i)} 
its $\in$-formalization can be expressed by a $\Delta_0$-formula;
\item  \label{cond:simplstprop(ii)} 
its meaning is invariant across forcing extensions: i.e. $P(\vec{a})$ holds in $V$ if and only if it holds in any generic extension of $V$;\footnote{By generic extension of $V$ we mean the standard way of coding a forcing extension of the universe of set within $V$ using, for example, a Boolean-valued construction.}
\item  \label{cond:simplstprop(iii)} 
the basic truths about $P(\vec{x})$ are forcing invariant, which roughly amounts to require that, whenever $\psi(\vec{x})$ is an $\in$-formula formalizing $P(\vec{x})$ and $\phi(\vec{y})$ is a boolean combination of $\Delta_0$-formulae and $\psi$, then the truth value of $\forall\vec{y}\,\phi(\vec{y})$ cannot be changed by forcing.\footnote{For example let $P(x)$ be the property \emph{$x$ is a stationary subset of $\omega_1$}, then $\exists x\, [P(x)\wedge P(\omega_1\setminus x)]$ is a $\ZFC$-theorem, hence its truth value cannot be changed by forcing. This is true regardless of the fact that whenever $G$ is $V$-generic for a forcing collapsing $\omega_1^V$ to become countable, the witnesses of the truth of $\exists x\, [P(x)\wedge P(\omega_1\setminus x)]$ in $V$ and that of $\exists x\, [P(x)\wedge P(\omega_1\setminus x)]$ in $V[G]$ cannot be the same set. Note that in an axiomatization of set theory in a signature $\tau$ which has stationarity on $\omega_1$ formalized by a unary predicate $P$ and has atomic predicates interpreting the $\Delta_0$-formulae, the above can be formalized by the $\Sigma_1$-sentence for $\tau$
 $\exists x\exists y\, [P(x)\wedge P(y) \wedge R(x,y)]$ with $R(x,y)$ being the atomic $\tau$-formula corresponding to the $\Delta_0$-property $(y=\cup x\setminus x)$.}
\end{enumerate}

Notice that the second and the third conditions are implied by the first. There are properties 
$P(x)$ which are not expressible using $\Delta_0$-formulae but for which nonetheless the second or third
(or both)
conditions apply; for example: the second condition applies to the provably $\Delta_1$-properties 
for $\ZFC^-$ and (assuming large cardinals)  to the property
\emph{``$x$ belongs to a certain universally Baire set $B$''}, while the third holds (again assuming large cardinals) for $P(x)$ being 
\emph{``$x$ is a stationary subset of $\omega_1$''} (see Thm. \ref{Thm: mainthmforcibility} and the comments following it).

We can give separate motivations for why each of the above condition is a simplicity criterion:
\ref{cond:simplstprop(i)} is motivated by set theoretic practice (we expand on this in Section \ref{subsec:rightsignst}).
\ref{cond:simplstprop(ii)} and \ref{cond:simplstprop(iii)} are
simplicity criteria as they state that we cannot use forcing to change the basic truths regarding properties satisfying them.
We consider forcing invariance a legitimate simplicity check, since forcing is the most efficient method to produce independence results and  (by the results we aim to expose in this paper) also the only method at all, at least for a large class of statements; we expand more on this in Section \ref{subsec:exclosedmodst}. Furthermore, notice that assuming large cardinals many set theoretic concepts become forcing invariant (this is the content of generic absoluteness results). In this respect if one is eager to accept large cardinal axioms as correct set theoretical principles (as we do), it is then natural to regard as ``simple'' all set theoretic concepts and properties which are made forcing invariant by large cardinals. 
\item \label{cond:maxstprop} 
Rather than focusing on the particularly delicate task of establishing which sets $Y$ do exist in $V$ and which do not (in order to argue whether $\exists x P(x,Y)$ can be consistent or valid), we just investigate if the $\Pi_2$-sentence $\forall x\exists y\,P(x,y)$ is consistent with the basic principles of set theory. This task can be realised by checking whether the sentence $\forall x\exists y\,P(x,y)$  belongs to $\bool{AMC}(T,A)$, for $T$ some $\in$-theory extending $\ZFC+$\emph{large cardinals}, and $A\subseteq\bool{Form}_{\bp{\in}}\times 2$ such that
$P(x,y)$ is formalizable by a quantifier free formula in signature $\bp{\in}_A$.
Moreover (assuming \ref{Descr:IWF}), if  $P(x,y)$ satisfies any of the above conditions for simplicity listed  in \ref{cond:simplstprop}, we can regard  a $\Pi_2$-sentence $\forall x\exists y P(x,y)$ which is consistent with the basic principles of set theory as a sentence expressing a true property of $V$. Indeed, by its being $\Pi_2$, such a sentence would assert that the (suitable fragment of the) universe of sets in which it holds is closed-off with respect to the
simple operations and relations described by $P$.
\end{enumerate}

Now that we have provided insights on which mathematical criteria we use to detect ``simplicity'' and that we have clarified  what we mean by ``consistent with the basic principles of set theory'', we can give a first general outline of the main results on the model companions of set theory.

\begin{enumerate}[(A)]
\item\label{cond:Amainres}
The (A)MC-spectrum of set theory can be used to characterize theories of the various $H_{\kappa^+}$ as $\kappa$ ranges over the infinite cardinals. More precisely, 
for any signature $\tau$ that extends $\bp{\in}$ and includes all $\Delta_0$-formulae among its atomic formulae, if the  $\tau$-theory of sets admits a model companion, this model companion is a model of 
Extensionality, Foundation, Union, it is closed under rudimentary functions or G\"odel operations, and satisfies the replacement axiom for a very large family of $\Sigma_1$-formulae for $\tau$. Furthermore, for any extensions $T$ of $\ZFC$ and any infinite cardinal $\kappa$ definable in $T$ we can (uniformly in $\kappa$) cook up at least one signature $\tau_\kappa$ 
(extending $\bp{\in}$ and including all $\Delta_0$-formulae among its atomic formulae) such that the $\tau_\kappa$-theory of $H_{\kappa^+}^{\mathcal{M}}$ common to all $\mathcal{M}$ which model $T$ is the AMC of $T$. Moreover these $\tau_\lambda$-theories of $H_{\lambda^+}$ cohere; more precisely for $\lambda<\kappa$ both definable in $T$ we can further assume that in any model $\mathcal{M}$ of $T$ the $\tau_\lambda$-theory of $H_{\lambda^+}^\mathcal{M}$ is correctly computed in the structure $(H_{\kappa^+}^\mathcal{M},\tau_\kappa^\mathcal{M})$, being $(H_{\lambda^+}^\mathcal{M},\tau_\lambda^\mathcal{M})$ a definable structure in $(H_{\kappa^+}^\mathcal{M},\tau_\kappa^\mathcal{M})$ (see Thm \ref{Thm:AMCsettheory+Repl}, and Notation \ref{not:modthnot2} for details).
\item\label{cond:Bmainres}
There is (more than one, and even recursive ones) $B\subseteq \bool{Form}_{\bp{\in}}\times 2$ such that:
\begin{itemize}
\item
$T_{\bp{\in},B}$ makes any $\Delta_0$-formula logically equivalent to some $\bp{\in}_B$-atomic formula, and
$\in_B$ has a constant symbol which is interpreted as $\omega_1^{\mathcal{M}}$ in any $\bp{\in}_B$-model $\mathcal{M}$ of $\ZFC+T_{\bp{\in},B}$.
\item
All quantifier free formulae of $\bp{\in}_B$ express simple properties according to simplicity criterion \ref{cond:simplstprop}
(assuming large cardinals). Furthermore the signature $\in_B$ introduces relations, constants, and function symbols which denote \emph{natural} and central concepts in the set-theoretical practice of the last fourty years.\footnote{The reader can find in \cite{VIABOURBAKI} a detailed introductory account of why the signatures $\in_B$ under consideration are \emph{natural}. A signature $\in_B$ which works can be obtained enriching the $\in$-signature with predicate symbols for the $\Delta_0$-formulae, function symbols for the G\"odel (and/or rudimentary) operations, constants for $\omega$ and $\omega_1$, a predicate symbol for the non-stationary ideal, and predicate symbols for the family of sets of reals definable in $L(\Ord^\omega)$ by an $\in$-formula without parameters; this family will give a recursive list of predicate symbols for universally Baire sets (assuming sufficiently many large cardinals exist). See also Def. \ref{not:STnotation} and/or \cite[Not. 3.1.1, Not. 5.1, Def. 5.2, Not. 5.3]{viale2021absolute}.}
\item
Whenever 
$S$ is any extension of the  $\bp{\in}$-theory
$\ZFC+$\emph{there are class many supercompact cardinals},\footnote{See Thm. \ref{Thm: mainthmforcibility} below for a precise formulation of the (as of now) optimal large cardinal assumptions.} $B$ in $\SpecAMC{S}{\bp{\in}}$.
\item
Whenever $S$ is as in the previous item,
the corresponding theory $\bool{AMC}(S,B)$ describes a theory of $H_{\omega_2}$, and is axiomatized by the $\Pi_2$-sentences (in the signature $\bp{\in}_B$) which hold in the $H_{\omega_2}$ of a model of strong forcing axioms.\footnote{E.g. strong forms of Woodin's axiom $(*)$ or $\MM^{++}$.} Among these $\Pi_2$-sentences there are $\neg\CH$ and a definable version $\theta$ of $2^{\aleph_0}=\aleph_2$.
\end{itemize}
\item \label{cond:Cmainres}
For\footnote{In this item \ref{cond:Cmainres} we can replace $\CH$ by an $\in$-sentence $\theta$ which is formalized by a $\Pi_2$-sentence  for $\in_B$ and is such that $\ZFC+\theta$ decides $2^{\aleph_0}=\aleph_2$; see Section \ref{mainthm:CH*} for details.} any theory $R$ extending $\ZFC$ such that $R+\neg\CH$
is consistent,
$\CH$  cannot be in $\bool{AMC}(R,A)$ for any $A$ in 
$\SpecAMC{R}{\bp{\in}}$ such that the $\Delta_0$-formulae are ($T_{\bp{\in},A}$-equivalent to) atomic $\bp{\in}_A$-formulae.
\item \label{cond:Dmainres}
For any\footnote{Item \ref{cond:Dmainres} applies to the $\Pi_2$-sentences $\neg\CH$ and $\theta$ considered in \ref{cond:Cmainres}.} $S,B$ as in \ref{cond:Bmainres}, any $\Pi_2$-sentence $\psi$ of 
$\bp{\in}_B$ is in $\bool{AMC}(S,B)$ if and only if $S$ proves that $\psi^{H_{\omega_2}}$ is forcible.
\item \label{cond:Emainres} 
The model completeness results for the theory of $H_{\aleph_2}$ of \ref{cond:Bmainres} mesh perfectly with certain forms of strong forcing axioms; for example one can prove that if one is in a model $(V,\in)$ of $S$ where the $\bp{\in}_B$-theory of $H_{\aleph_2}^V$ is model complete (hence the AMC of the $\in_B$-theory of $V$), then a strong form of bounded Martin's maximum holds.
\end{enumerate}

%

The resulting picture of $V$ that emerges from \ref{cond:Amainres} is that of a cumulative hierarchy whose initial segments (the  $H_{\kappa^+}$'s) are closed off with respect to more and more complex set-theoretic operations and produce a coherent step by step approximations of the theory of $V$ given by the various theories of $H_{\lambda^+}$ as $\lambda$ ranges over the infinite cardinals.

In this way, set theory is treated algebraically like field theory: once we choose the cardinal $\kappa$ and a corresponding signature $\tau_\kappa$, the initial segment $H_{\kappa^+}$ of the universe $V$ has a theory which behaves, with respect to $\ZFC_{\tau_\kappa}$ (e.g. $\ZFC$ enriched with the replacement axiom for $\tau_\kappa$-formulae, see \ref{not:STnotation} for details), as algebraically closed fields behave with respect to field theory as axiomatized in $\bp{\cdot,+,0,1}$. 

Moreover, by \ref{cond:Bmainres} above, this analogy between set theory and algebra is exemplified for the cardinal $\aleph_1$ by a signature $\in_B$ which renders the theory of $H_{\aleph_2}$ the result of a closing-off process with respect to simple set-theoretical concepts (according to  our simplicity criteria \ref{cond:simplstprop}), and to the standard set theoretic practice of the last fourty years.\footnote{Universally Baire sets and the non-stationary ideal are regarded as tame concepts (at least in the presence of large cardinals), and $\Delta_0$-properties are certainly considered elementary set-theoretic concepts, see Section \ref{subsec:rightsignst} below.} Consequently, the algebraic completeness displayed by $H_{\aleph_2}$, as the model companion of set theory with respect to $\in_B$, satisfies the conditions described in \ref{Descr:IWF} and substantiates the Hilbertian completeness encoded by it. Furthermore, \ref{cond:Cmainres} shows that the study of this  model companion of set theory gives trustworthy information on the value of the continuum, while \ref{cond:Dmainres} suggests an explanation of the incredible success of the method of forcing (and of forcing axioms) in providing a consistent solution to problems formalizable at the level of $H_{\aleph_2}$ as $\Pi_2$-sentences for $\in_B$. Finally condition  \ref{cond:Emainres} shows that the present model-theoretic approach and the use of forcing axioms are two sides of the same coin. This, in turn, can be seen as a further justification for extending the axioms of set theory with forcing axioms and, consequently, a justification of G\"odel's program at least at the level of (a large fragment of) third order arithmetic.

\subsection{What is the right signature for set theory?}\label{subsec:rightsignst}

The standard axioms of $\ZFC$ in the $\in$-signature are clearly sufficient to provide a first order axiomatization of set theory. However a closer inspection reveals that many simple set-theoretic concepts are not formalized by simple $\in$-formulae. 

Consider, for example, the notion of ordered pair. While we informally write $x=\ap{y,z}$ to mean that
\emph{$x$ is the ordered pair with first component $y$ and second component $z$}, in set theoretic terms this statement hides a non-trivial coding of the concept of ordered pair (for example) by means of Kuratowski's definition: $x=\bp{\bp{y},\bp{y,z}}$. A proper definition of the concept of ordered pair in the $\in$-signature can then be given by the following $\in$-formula:
\[
\exists t\exists u\;[\forall w\,(w\in x\leftrightarrow w=t\vee w=u)
\wedge\forall v\,(v\in t\leftrightarrow v=y)
\wedge\forall v\,(v\in u\leftrightarrow v=y\vee v=z)].
\]

It is clear that the meaning of this $\in$-formula is hardly recognizable with a rapid glance (unlike $x=\ap{y,z}$). Moreover, from a purely logical perspective, its L\'evy complexity  is already $\Sigma_2$. This clashes with our understanding that the concept of ordered pair is simple. Indeed, we do not regard the notion of ordered pair as a complex concept, contrary to other more complicated and theoretically loaded ones like that of uncountability, or many of the properties of the continuum (such as its correct place in the hierarchy of uncountable cardinals). In a similar vein other very basic notions such as being a function, a binary relation, or the domain or the range of a function are formalized by rather complicated $\in$-formulae, both from the point of view of their readability and of their L\'evy complexity.

The standard solution adopted in set theory textbooks\footnote{See for example \cite[Chapter IV, Def. 3.5]{KUNEN} or \cite[Def. 12.9]{JECHST}.} is to regard as basic all
those $\in$-formulae in which the quantifiers are bounded to range
over the elements of some set, that is the $\Delta_0$-formulae. 
In order to make these observations precise we need to be extremely cautious on our notational conventions.

\begin{notation} \label{not:STnotation}
For any 
$A\subseteq \bool{Form}_{\bp{\in}}\times 2$ 
we write $\in_A$ rather than $\bp{\in}_A$, and we let $T_{\in,A}$ be the $\in_A$-theory
\[
T_{\bp{\in},A}+\forall x \,\qp{(\forall y\,y\notin x)\leftrightarrow c_{\bp{\in}}=x},
\]
where the theory $T_{\bp{\in},A}$ (according to Notation \ref{not:keynotation0} for $\bp{\in}$ and $A$) is reinforced  by  
an axiom asserting that the interpretation of the constant symbol $c_{\bp{\in}}$ is the empty set.
We use the abbreviations $\in_{\Delta_0}$ and $T_{\Delta_0}$ to denote what, according to the above conventions, should rather be slight extensions of\footnote{The precise definition of $\in_{\Delta_0}$ and $T_{\Delta_0}$ can be found in \cite[Notation 3.1.1]{viale2021absolute}. They are obtained by enriching $\in_{\Delta_0\times \bp{0}}$ and $T_{\in,\Delta_0\times \bp{0}}$ with  a constant symbol for $\omega$ and function symbols for the G\"odel operations (as defined in \cite[Def. 13.6]{JECHST}) and axioms to interpret them correctly.} 
$\in_{\Delta_0\times \bp{0}}$ and $T_{\in,\Delta_0\times \bp{0}}$.

For any $\tau\supseteq\in_{\Delta_0}$ we also let
$\ZFC_\tau$ denote $\ZFC+T_{\Delta_0}$ enriched with the replacement axiom\footnote{We use the following strong form of the replacement axiom:
whenever $R\subseteq X\times V$ is a definable class with domain $X$ being a set, there is a uniformizing function $f:X\to V$ with $f$ a set and $R(x,f(x))$ holding for all $x\in X$. This strong form of replacement holds in $\in$-models $(V,\in)$ of $\ZFC$ and is inherited by the substructures of the form $(H_\lambda,\in)$ for $\lambda$ a regular cardinal of $V$. 
The resulting weak set theory one obtains when dropping the powerset axiom does not have the pathological models which can have other formalizations of set theory without powerset axiom where replacement is expressed by a weaker axiom schema (see for example \cite{hamkins:zfc_p}).} for all $\tau$-formulae, and $\ZFC^-_\tau$ denote
$\ZFC_\tau$ without the powerset axiom. $\ZFC_{\Delta_0}$ denotes $\ZFC_\tau$ for 
$\tau$ being $\in_{\Delta_0}$, and similarly we define $\ZFC^-_{\Delta_0}$.
\end{notation}

The reason why set-theoretic practice regards
the concepts expressed by quantifier free formulae of $\in_{\Delta_0}$ as ``simple'' is the fact that the truth value of these formulae is invariant among transitive models of large enough fragments of $\ZF$ (e.g. \cite[Corollary IV.3.6]{KUNEN}), and thus also forcing invariant (e.g. \cite[Lemma 14.21]{JECHST}). Furthermore $\in_{\Delta_0}$ is a signature which allows one  to formalize many fundamental set theoretic concepts using formulae whose L\'evy complexity is in accordance to our intuitive understanding 
(e.g. \cite[Chapter 13, Lemma 13.10]{JECHST}).

For example consider the following notions and their logical complexity:

\begin{itemize}
\item
$(x\text{ is a cardinal})$ is the $\Pi_1$-formula (for $\in_{\Delta_0}$)
\[
(x\text{ is an ordinal})\wedge
\forall f\,\qp{(f\text{ is a function}\wedge \dom(f)\in x)\rightarrow\ran(f)\neq x}.
\]
\item
$(x\emph{ is }\aleph_1)$ is the boolean combination of $\Sigma_1$-formulae
\begin{align*}
(x\text{ is a cardinal})\wedge (\omega\in x) \wedge\\
\wedge\exists F\,\qp{(F:\omega\times x\to x)\wedge \forall\alpha\in x\,(F\restriction \omega\times\bp{\alpha}\text{ is a surjection on }\alpha)}.
\end{align*}
\item
$\CH$ is the $\Sigma_2$-sentence
\[
\exists f\, \qp{(f\text{ is a function}\wedge\dom(f) \emph{ is }\aleph_1) \wedge \forall r\subseteq\omega\, (r\in\ran(f))}.
\]
and $\neg\CH$ is the  boolean combination of $\Pi_2$-sentences\footnote{We let $\neg\CH$ include the $\Sigma_2$-sentence $\exists x\,(x\emph{ is }\aleph_1)$, for otherwise its failure could be witnessed by the assertion that there is no uncountable cardinal, a statement which holds true in $H_{\omega_1}$, regardless of whether $\CH$ or its negation is true in the corresponding universe of sets.}
\[
\exists x\,(x\emph{ is }\aleph_1)\wedge\forall f\, \qp{(\dom(f)\emph{ is }\aleph_1\wedge f\text{ is a function}) \rightarrow \exists r\subseteq\omega\, (r\not\in\ran(f))}.
\]
\item
$(x\emph{ is }\aleph_2)$ is the $\Sigma_2$-formula
\begin{align*}
(x\text{ is a cardinal})\wedge \\
\wedge\exists F\exists y\,\qp{ (y\emph{ is }\aleph_1)\wedge (y\in x)\wedge(F:y\times x\to x)\wedge \forall\alpha\in x\,(F\restriction y\times\bp{\alpha}\text{ is a surjection on }\alpha)}.
\end{align*}
\item
$2^{\aleph_0}>\aleph_2$ is the boolean combination of $\Pi_2$-sentences
\[
\exists x\,(x\emph{ is }\aleph_2)\wedge\forall f\, \qp{(f\text{ is a function}\wedge\dom(f)\emph{ is }\aleph_2) \rightarrow \exists r\,(r\subseteq\omega \wedge r\not\in\ran(f))}.
\]
\item
$2^{\aleph_0}\leq\aleph_2$ is the $\Sigma_2$-sentence
\[
\exists f\, \qp{(f\text{ is a function})\wedge\dom(f)\emph{ is }\aleph_2 
\wedge \forall r\, (r\subseteq\omega\rightarrow r\in\ran(f))}.
\]
\end{itemize}

Let us also introduce the notation we will use to handle the substructure relation over expanded signatures.
 The following conventions supplement Notation \ref{not:keynotation0}. 
 \begin{notation}
 Given $\tau$-structures $\mathcal{M}$, $\mathcal{N}$,
$\mathcal{M}\sqsubseteq \mathcal{N}$ denotes the substructure relation,
$\mathcal{M}\prec_n \mathcal{N}$ denotes the $\Sigma_n$-substructure relation, and 
$\mathcal{M}\prec \mathcal{N}$ denotes the elementary substructure relation.
 \end{notation}
 
\begin{notation}\label{not:modthnot2}
Let $\tau\supseteq\in_{\Delta_0}\cup\bp{\kappa}$ be a signature with $\kappa$ a constant symbol, 
$(M,\tau^M)$ a
$\tau$-structure, and  $B$  a subset of $\bool{Form}_{\tau}\times 2$. Then,
$(M,\tau_B^M)$ is the unique extension of 
$(M,\tau)$ defined in accordance with 
Notations \ref{not:keynotation0} 
which satisfies $T_{\tau,B}$.
In particular $(M,\tau_B^M)$ is a shorthand for 
$(M,S^M:S\in\tau_B)$.
If $(N,\tau^N)$ is a substructure of $(M,\tau^M)$ we also write
$(N,\tau_B^M)$ as a shorthand for 
$(N,S^M\restriction N:S\in\tau_B)$.

\end{notation}

\subsection{L\'evy absoluteness and the possible AMCs of set theory}\label{subsec:levabsAMCST}

Recall that for an infinite cardinal $\lambda$, $H_{\lambda}$ is the initial transitive fragment of $V$ given by those sets whose transitive closure has size less than $\lambda$ and is by itself a set in $V$
(see \cite[Section IV.6]{KUNEN}).
We have the following results.

\begin{description} 
\item[Stratification of $V$] the universe of sets $V$ can be stratified as the union, along the class of infinite cardinals $\kappa$, of the sets $H_{\kappa^+}$.
\item[Standard theory of $H_{\lambda}$] \cite[Thm. IV.6.5]{KUNEN} for any uncountable cardinal $\lambda$, $(H_{\lambda},\in_{\Delta_0}^V)$ is a model of $\ZFC^-_{\Delta_0}$ (recall Notation \ref{not:STnotation}).
\item[Strong L\'evy absoluteness]  \cite[Lemma 5.3]{VIAVENMODCOMP} For all $A_i\subseteq \pow{\kappa}^{n_i}$ for $i=1,\dots,k$,
\[
(H_{\kappa^+}^V,\in_{\Delta_0}^V,A_1,\dots,A_k)\prec_1(V,\in_{\Delta_0}^V,A_1,\dots,A_k).
\] 
\item[Second order characterization of $H_{\kappa^+}$] Whenever $\mathcal{M}$ is an $\in_{\Delta_0}\cup\bp{\kappa}$-model of $\ZFC_{\Delta_0}+\kappa$\emph{ is an infinite cardinal}
$(H_{\kappa^+}^{\mathcal{M}},\in_{\Delta_0}^{\mathcal{M}},\kappa^{\mathcal{M}})$ is a model of the $\Pi_2$-sentence for 
$\in_{\Delta_0}\cup\bp{\kappa}$
\begin{equation}\label{eqn:charHkappa+}
\forall x\exists f\,(f:\kappa\to x\text{ is a surjection}).
\end{equation}
Furthermore $H_{\kappa^+}^V$ can be described as the unique set $M\in V$ such that;
\begin{itemize}
\item
$\kappa\in M$ and $M$ is transitive;
\item
$(M,\in_{\Delta_0}^V)$ satisfies (\ref{eqn:charHkappa+}) and $\ZFC^-_{\Delta_0}$;
\item
$M$ has the (second order property) that for all $a\in M$ and $b\subseteq a$, $b\in M$ as well.
\end{itemize}
\end{description}

In particular the theory of each $H_{\kappa^+}$ offers a $\Sigma_1$-approximation of the theory of $V$ with respect to a signatures $\tau_\kappa$ which extends
$\in_{\Delta_0}\cup\bp{\kappa}$ and contains predicates for subsets of $\pow{\kappa}^n$ for some $n\in\mathbb{N}$. Consequently, if we want to study the $\Sigma_1$-properties of the reals (second order arithmetic), or the powerset of the reals (third order arithmetic) it is sufficient to study them within the theory of some  $H_{\kappa^+}$ for $\kappa$ a large enough regular cardinal (and for most purposes $H_{\aleph_1}$ suffices for second order arithmetic, and $H_{\aleph_2}$ for third-order arithmetic). 
Moreover, by L\'evy absoluteness, $H_{\kappa^+}$ and $V$ agree on the computation of many simple set theoretic operations and relations: for example the operations and relations formalizable by means of $\Delta_0$-formulae and $\Delta_0$-definable Skolem functions.

The above results entail that whenever $T$ is a $\tau$-theory that extends  $\ZFC_{\Delta_0}+\kappa$\emph{ is a cardinal},   $\tau$ is a signature that contains $\in_{\Delta_0}\cup\bp{\kappa}$
and is contained in 
\[
\in_{\Delta_0}\cup\bp{A: A\subseteq\pow{\kappa}^n,\, n\in\mathbb{N}},
\] 
and $(V,\tau^V)$ models $T$, then  
\[
(H_{\kappa^+}^V,\tau^V)\prec_1 (V,\tau^V).
\]
Hence
 $(H_{\kappa^+},\tau^V)$ witnesses that (\ref{eqn:charHkappa+})
is consistent with the universal fragment of the $\tau$-theory of $V$.
In particular if $T$ has an AMC for $\tau$, say $S$, then  (\ref{eqn:charHkappa+})
must be in $S$, since it is a $\Pi_2$-sentence for $\tau$ which is strongly $T^\tau_{\forall\vee\exists}$-consistent.


We can now collect all these observations and conclude the following.
\begin{enumerate}
\item The structure of $V$ is uniquely determined by the structure of the various $H_{\kappa^+}$ as $\kappa$ ranges among the infinite cardinals.
\item Given an infinite cardinal $\kappa$, consider a signature $\tau$ extending $\in_{\Delta_0}$ only with predicates
for subsets of $\pow{\kappa}$ and a constant symbol for $\kappa$. Then, if some $T\supseteq\ZFC_\tau$
has an AMC, this AMC looks like a theory of $H_{\kappa^+}$ and is axiomatized by the $\Pi_2$-sentences for $\tau$ which are strongly  $T^\tau_{\forall\vee\exists}$-consistent.
\end{enumerate}

Now assume $T$ is some $\in$-theory $T\supseteq\ZFC+$\emph{large cardinals} and $A\subseteq \bool{Form}_\in\times 2$ is in $\SpecAMC{T}{\bp{\in}}$.
Assume further that:
\begin{itemize}
\item
 for some constant symbol $\kappa$ of $\in_A$
\[
\ZFC^-+\forall X\exists f\,(f:\kappa\to X\text{ is a surjection})
\]
is in $\bool{AMC}(T,A)$;
\item 
every atomic formula of
$\in_A$ satisfies at least one of the simplicity criteria set forth in 
\ref{cond:simplstprop}. 
\end{itemize}
Then (on the basis of \ref{Descr:IWF}) we should conclude that $\bool{AMC}(T,A)$ describes the correct theory of $H_{\kappa^+}$, since criterion \ref{cond:maxstprop} is satisfied by all axioms of $\bool{AMC}(T,A)$.


\subsection{Existentially closed fragments of the set theoretic universe versus AMC} \label{subsec:exclosedmodst}

We conclude this section by giving a non-exhaustive list of results showing that, for certain infinite cardinals $\kappa$, the corresponding $H_{\kappa^+}$ already witnesses some important existential closure properties (e.g. L\'evy absoluteness, Shoenfield's absoluteness, $\mathsf{BMM}$, $\mathsf{BMM}^{++}$) and its theory has some traits which are peculiar of model complete theories         (e.g. Woodin's absoluteness, $\MM^{+++}$).
 The reader does not need to be familiar with these results, as they only serve as motivation to bring our focus on existentially closed models for set theory.\footnote{For further examples of existentially closed models of set theory see \cite{Venturi2020-VENIFA}.}



\begin{description}
\item[L\'evy absoluteness] \cite[Lemma 5.3]{VIAVENMODCOMP}
Whenever $\lambda$ is a regular uncountable cardinal,
\[
(H_\lambda,\in_{\Delta_0}^V)\prec_1 (V,\in_{\Delta_0}^V).
\] 
 \item[Shoenfield's absoluteness]  (see \cite[Lemma 1.2]{VIAMMREV} for the apparently weaker formulation we give here)
  Whenever $G$ is $V$-generic for some forcing notion in $V$,
 \[
 (H_{\omega_1},\in_{\Delta_0}^V)\prec_1 (V[G],\in_{\Delta_0}^{V[G]}).
 \]
  \item[Woodin's absoluteness]  (see \cite[Lemma 3.2]{VIAMMREV} for the weak form of Woodin's result we give here) 
  Whenever $G$ is $V$-generic for some forcing notion in $V$ (and there are class many Woodin cardinals in $V$),
  \[
  (H_{\omega_1}^V,\in_{\Delta_0}^V)\prec (H_{\omega_1}^{V[G]},\in_{\Delta_0}^{V[G]}).
  \]
   \item[Bounded Martin's Maximum ($\mathsf{BMM}$)]  \cite{BAG00}
   Whenever $G$ is $V$-generic for some stationary set preserving forcing notion in $V$,
   \[
   (H_{\omega_2},\in_{\Delta_0}^V)\prec_1 (V[G],\in_{\Delta_0}^{V[G]}).
   \]
     \item[$\mathsf{BMM}^{++}$]   \cite[Def. 10.91]{WOOBOOK}
     Whenever $G$ is $V$-generic for some stationary set preserving forcing notion in $V$,
      \[
      (H_{\omega_2},\in_{\Delta_0}^V,\NS_{\omega_1}^V)
      \prec_1 (V[G],\in_{\Delta_0}^{V[G]},\NS_{\omega_1}^{V[G]}),
      \]
       where $\NS_{\omega_1}$ is a unary predicate symbol interpreted by the 
       non-stationary ideal on $\omega_1$.
      \item[Bounded category forcing axioms, $\MM^{+++}$, $\bool{RA}_\omega(\SSP)$]
       \cite{VIAASP,VIAAUD14,VIAMM+++}
       Whenever 
      $V$ and $V[G]$ are models of $\MM^{+++}$ ($\bool{RA}_\omega(\SSP)$, $\bool{BCFA}(\SSP)$) 
      and $G$ is $V$-generic for some stationary
	set preserving forcing notion in $V$,
      \[
      (H_{\omega_2}^V,\in_{\Delta_0}^V)\prec (H_{\omega_2}^{V[G]},\in_{\Delta_0}^{V[G]}).
      \]
      \item[Woodin's axiom $(*)$] (see \cite[Def. 7.5]{HSTLARSON}, or \cite{ASPSCH(*)} for the formulation given here) 
      Whenever there are class many Woodin cardinals, $\NS_{\omega_1}$ is a precipitous ideal, and $G$ is $V$-generic for some stationary set preserving forcing notion in $V$
\[
(H_{\omega_2},\in_{\Delta_0},\NS_{\omega_1}, A: A\in \pow{\pow{\omega}}^{L(\mathbb{R})})\prec_1  (V[G],\in_{\Delta_0}^{V[G]},\NS_{\omega_1}^{V[G]}, A^{V[G]}: A\in \pow{\pow{\omega}}^{L(\mathbb{R})})
\]
(where $A^{V[G]}$ is the canonical interpretation in $L(\mathbb{R})^{V[G]}$ of the definition of $A$ in $L(\mathbb{R})$).
     \end{description}


   \section{Main results} \label{sec:mainresults}
We can now present results culminating  the discussion of the previous sections. The proofs can be found in \cite{viale2021absolute}.
   
   \begin{notation}
%

For $A\subseteq \bool{Form}_\in\times 2$, $\ZFC_A$ denotes the $\in_A$-theory 
$\ZFC_{\in_A}$ defined in Notation \ref{not:STnotation}. Similarly we define $\ZFC^-_A$.
\end{notation}

\begin{definition}
Let $T\supseteq\ZFC^-$ be an $\in$-theory. $\kappa$ is a $T$-definable cardinal if for some $\in$-formula $\phi_\kappa(x)$, $T$ proves:
\begin{itemize} 
\item
$\exists!x\,\phi_\kappa(x)$ and
\[
\forall x\,[\phi_\kappa(x)\rightarrow (x\emph{ is a cardinal})].
\]
\item
$\kappa$ is the constant $f_{\phi_\kappa}$ existing in the signature $\in_{\bp{(\phi_\kappa,1)}}$.
\end{itemize} 
\end{definition}

\subsection{The AMCs of set theory are theories of $H_{\kappa^+}$}
The first result shows that the AMC spectrum of set theory isolates a rich set of theories which produce models of
$\ZFC^-$, that is, structures which behave like an $H_{\lambda}$ for some regular uncountable cardinal $\lambda$.

\begin{theorem}\label{Thm:AMCsettheory+Repl}
Let $R$ be an $\in$-theory extending $\ZFC$. 
\begin{enumerate}[(i)]
\item \label{Thm:AMCsettheory+Repl-1}
Assume $A\in \SpecMC{R}{\in}$ and $\in_A\supseteq \in_{\Delta_0}$. Then any model of $\bool{MC}(R,A)$ models extensionality, foundation, choice, is closed under G\"odel operations, and satisfies the replacement axiom for any $\Sigma_1$-formula $\psi(x,y)$ of $\in_{\Delta_0}$ such that $R+T_{\in,\A}$ models
$\forall x\exists y\psi(x,y)$.
\item
Assume $A\in \SpecMC{R}{\in}$, $\in_A\supseteq \in_{\Delta_0}$, and $\lambda$ is an $R$-definable cardinal represented by a constant symbol of $\in_A$ and such that\footnote{$H_{\lambda^+}^\mathcal{M}$ denotes the substructure of $\mathcal{M}$ whose extension is given by the formula defining $H_{\lambda^+}$ in the model (using the parameter $\lambda$).}
\[
(H_{\lambda^+}^{\mathcal{M}},\in_A^{\mathcal{M}})\prec_1\mathcal{M}
\] 
whenever $\mathcal{M}$ models $R+T_{\in,A}$. Then $\forall x\exists f\,(f:\lambda\to x\text{ is a surjection})$ is in $\bool{MC}(R,A)$.
\item \label{Thm:AMCsettheory+Repl-2}
Assume $\lambda$ is an $R$-definable cardinal.\footnote{Furthermore we can uniformly choose the sets $A_\lambda$ so that in any model $\mathcal{M}$ of $R$ and for $\kappa,\lambda$ definable $R$-cardinals such that $\mathcal{M}$ models $\kappa>\lambda$,  $(H_{\lambda^+}^{\mathcal{M}},\in_{A_\lambda}^{\mathcal{M}})$ is a definable structure in $(H_{\kappa^+}^{\mathcal{M}},\in_{A_\kappa}^{\mathcal{M}})$ (more precisely, for each finite $\sigma\subseteq \in_{A_\lambda}$,  $(H_{\lambda^+}^{\mathcal{M}},\sigma^{\mathcal{M}})$ is definable in $(H_{\kappa^+}^{\mathcal{M}},\in_{A_\kappa}^{\mathcal{M}})$).}
Then there exists $A_\lambda\in \SpecAMC{R}{\in}$ with $\in_{A_\lambda}$ containing $\in_{\Delta_0}$ such that 
$\bool{AMC}(R,A_\lambda)$ is given by the $\in_{A_\lambda}$-theory common to the structures $H_{\lambda^+}^\mathcal{M}$ as $\mathcal{M}$ ranges among the $\in_{A_\lambda}$-models of $R+T_{\in,A_\lambda}$.
\end{enumerate}
\end{theorem}

\subsection{Forcibility versus absolute model companionship}
The following is the major result relating AMC to forcibility and to  forcing axioms:\footnote{The reader unaware of what $\MM^{++}$ or a stationary set preserving forcing is can skip the second and third items of the theorem.}

\begin{theorem}\label{Thm: mainthmforcibility}
Let $S$ be the $\in$-theory 
\[
\ZFC+\emph{there exist class many supercompact cardinals}.
\]
Then there is a recursive set\footnote{$\in_B$ can be for example $\in_{\NS_{\omega_1},\mathcal{A}}$  as defined in \cite[Notation 5.3]{viale2021absolute} where $\mathcal{A}$ consists of the sets of reals lightface definable  in $L(\Ord^\omega)$; see also \ref{cond:Bmainres} on p. \pageref{cond:Bmainres}.} $B\in\SpecAMC{S}{\in}$ with $\in_B$ containing $\in_{\Delta_0}$, and such that for any 
$\Pi_2$-sentence $\psi$ for $\in_B$ and any $\in$-theory $R\supseteq S$ the following are equivalent:
\begin{enumerate}[(a)]
\item \label{Thm: mainthmforcibility-2-a} $\psi\in \bool{AMC}(R,B)$;
\item \label{Thm: mainthmforcibility-2-b}
$(R+T_{\in,B})^{\in_B}_{\forall\vee\exists}+S+\MM^{++}+T_{\in,B}$ proves\footnote{Here and elsewhere we write $\psi^N$ to denote the relativization of $\psi$ to a definable class (or set) $N$; see \cite[Def. IV.2.1]{KUNEN} for details. Recall that for a $\tau$-theory $S$ $S^\tau_{\forall\vee\exists}$ denotes the boolean combinations of universal $\tau$-sentences which follow from $S$.} $\psi^{H_{\omega_2}}$;
\item 
$R$ proves that  $\psi^{H_{\omega_2}}$ is forcible\footnote{Here and in the next item we mean that the $\in$-formula $\theta$ which is $T_{\in,B}$-equivalent to $\psi$ is such that $\theta^{H_{\omega_2}}$ is forcible by the appropriate forcing.} by a stationary set preserving forcing;
\item\label{Thm: mainthmforcibility-2-c}
$R$ proves that  $\psi^{H_{\omega_2}}$ is forcible by some forcing;
\item For any $R'\supseteq R$, $\psi+(R'+T_{\in,B})^{\in_B}_{\forall\vee\exists}$ is consistent.
\end{enumerate}

Furthermore for any $\theta$ which is a boolean combination of $\Pi_1$-sentences for $\in_B$ and any
$(V,\in)$ model of $S$, TFAE:
\begin{enumerate}[(A)]
\item \label{Thm: mainthmforcibility-2-A}
$(V,\in_B^V)$ models $\theta$;
\item \label{Thm: mainthmforcibility-2-B}
$(V,\in_B^V)$ models that some forcing notion $P$ forces $\theta$;
\item \label{Thm: mainthmforcibility-2-C}
$(V,\in_B^V)$ models that all forcing notions $P$ force $\theta$.
\end{enumerate}
\end{theorem}
The equivalence of \ref{Thm: mainthmforcibility-2-a} and \ref{Thm: mainthmforcibility-2-b} shows that the $\in_B$-theory of $H_{\aleph_2}$ is model complete in models of strong forcing axioms; we give a partial converse in Thm. \ref{thm:modcompBMM}.

The second part of the theorem shows that forcing cannot change the $\Pi_1$-fragment of the theory of $V$ in signature ${\in_B}\supseteq{\in_{\Delta_0}}$.
Moreover, a direct consequence of the above result is that if $(V,\in)$ is a model of $S$ and $S^*$ is the $\in_B$-theory of its unique extension to a model of 
$T_{\in,B}$, we get that a $\Pi_2$-sentence $\psi$ for $\in_B$ is consistent with the universal fragment of $S^*$
if and only if $\psi^{H_{\omega_2}}$ is forcible over $V$.

\subsection{The AMC-spectrum of set theory and the continuum problem}\label{subsec:AMCCH}

The study of the AMC spectrum of set theory places $\aleph_2$ in a very special position among the possible values of the continuum.

\begin{theorem} \label{mainthm:CH*}
The following holds:
\begin{enumerate} 
\item \label{mainthm:CH*-1}
Let  $R\supseteq \ZFC$ be an $\in$-theory. 
Assume $A\in \SpecAMC{R}{\in}$ is such that $\in_A$ contains $\in_{\Delta_0}$ and
$\neg\CH+(R+T_{\in,A})^{\in_A}_{\forall\vee\exists}$ is consistent. 
Then $\CH\not\in \bool{AMC}(R,A)$.
\item \label{mainthm:CH*-2}
For the signature $\in_B$ and the the theory $S$ of Thm. \ref{Thm: mainthmforcibility}
$\neg\CH$ is in $\bool{AMC}(R,B)$ for any $\in$-theory $R\supseteq S$.
\end{enumerate}
\end{theorem}

Let us briefly argue why $\CH$ cannot be in $\bool{AMC}(R,A)$ whenever $R+\neg\CH$ is consistent. 
\begin{proof}
We use the following peculiar property of AMC (which can fail for model companionship\footnote{For example no algebraically closed field can be a model of the $\Pi_2$-sentence (actually $\Pi_1$) $\forall x\,(\neg x^2+1=0)$ (which is consistent with the theory of fields, as it holds in $\mathbb{Q}$).}):
\begin{fact}\cite[Proposition 2.2.8, Def. 2.2.4, Lemma 2.3.4]{viale2021absolute} 
Let $T,S$ be $\tau$-theories with $T$ the AMC of $S$, and $\psi$ be a $\Pi_2$-sentence such that $S+\psi$ is consistent. Then there is a model of $T+\psi$.
\end{fact}


Recall that $\CH$ is the conjunction of two $\Sigma_2$-sentences (for $\in_{\Delta_0}$): $\psi_0$ stating ``\emph{the first uncountable cardinal exists}'' and $\psi_1$ stating ``\emph{there is a surjective function with domain the first uncountable cardinal and range  the reals}'', while $\neg\CH$ is the conjunction of $\psi_0$ with the $\Pi_2$-sentence $\neg\psi_1$. 

Now if $\psi_0\notin\bool{AMC}(R,A)$, $\CH\not\in\bool{AMC}(R,A)$ as well. Otherwise
$R+\neg\CH$ is consistent by assumption,
hence so is 
$R+T_{\in,A}+\psi_0\wedge\neg\psi_1$, since the latter is a conservative extension of the former.
By the above, we can find a model $\mathcal{M}$ of $\bool{AMC}(R,A)+\neg\psi_1$. This model witnesses 
that $\bool{AMC}(R,A)$ cannot prove $\CH$, even when $\psi_0\in\bool{AMC}(R,A)$.
\end{proof}

Let us also argue why
$\neg\CH$ falls in the AMC of $R+T_{\in,B}$ for the set $B\subseteq\bool{Form}_\in\times 2$ and the theory $R$ of (\ref{mainthm:CH*-2}) above. Towards this aim we appeal to the equivalence between \ref{Thm: mainthmforcibility-2-A} and \ref{Thm: mainthmforcibility-2-C} of Thm. \ref{Thm: mainthmforcibility} (which gives a precise instantiation of simplicity criterion \ref{cond:simplstprop(iii)}).
First of all notice  that an $\in$-theory $R$ is complete if and only if so is $R+T_{\in,B}$.
Let $R_0$ be some $\in$-completion of $R\supseteq S$ and $(V,\in)$ a model of $R_0$. 
We now note that $\in_B$ has a constant $\kappa$ to denote $\omega_1$, $\bool{AMC}(S,B)$ satisfies the $\in_{\Delta_0}\cup\bp{\kappa}$-sentence \emph{``$\kappa$ is the first uncountable cardinal''}, $\bool{AMC}(R_0,B)=\bool{AMC}(S,B)+(R_0+T_{\in,B})_{\forall\vee\exists}$ (by Lemma \ref{lem:charAMCvsMC}).
Let $\bool{B}$ be a cba such that $\Qp{\neg\CH}_{\bool{B}}=1_{\bool{B}}$ holds in $(V,\in)$.\footnote{For example the boolean completion of the partial order originally devised by Cohen. See for details \cite[Thm. 14.32]{JECHST}}
Let $G$ be any ultrafilter on $\bool{B}$. Then 
\[
V^{\bool{B}}/_G\models\forall f (f \text{ is a function with domain }\kappa\rightarrow \exists r \subseteq\omega \text{ which is not in its range})
\]
(by \cite[Lemma 14.14]{JECHST}) 
and $V^{\bool{B}}/_G$ also models the universal and existential $\in_B$-fragment of $R+T_{\in,B}$ (by 
\ref{Thm: mainthmforcibility-2-A}$\Leftrightarrow$\ref{Thm: mainthmforcibility-2-C} of Thm. \ref{Thm: mainthmforcibility}, and \cite[Lemma 14.14]{JECHST}).
In particular, the $\Pi_2$-conjunct of $\neg\CH$ is consistent with the universal and existential fragment of any completion of $S+T_{\in,B}$, hence it belongs to the AMC of $S+T_{\in,B}$, while the $\Sigma_2$-conjunct of $\neg\CH$ is in the AMC of $S+T_{\in,B}$ since the axiom $\bool{Ax}^1_\psi$ for $\psi$ defining the first uncountable cardinal is in $T_{\in,B}$.

\medskip

In particular we see that forcing becomes a powerful tool to prove that a $\Pi_2$-sentence formalizable in
$\in_B$ is in the AMC of $S+T_{\in,B}$. Indeed, it suffices to prove over $S$ that this sentence is forcible.
More generally minimal variants of the above argument give:

\begin{remark}
Let  $S$ be the theory of Thm. \ref{Thm: mainthmforcibility}, and $R_0\supseteq\ZFC$ be such that $R_0\cup S$ is consistent. Consider an $\in$-sentence $\psi$ such that $\psi$ as formulated in signature $\in_B$ is $R_0+T_{\in_B}$-equivalent to
a $\Pi_2$-sentence for $\in_B$ such that $R_0$ proves that $\psi^{H_{\aleph_2}}$ is forcible.
Then:
\begin{itemize} 
\item
$\psi\in\bool{AMC}(R,B)$ for any $R\supseteq S\cup R_0$;
\item 
$\neg\psi\not\in\bool{AMC}(T,A)$ for any $A$ and $\in$-theory $T\supseteq R_0$ such that $\in_A$ contains all the symbols of $\in_B$ appearing in $\psi$, 
and $(T+T_{\in_A})_{\forall\vee\exists}+\psi$ is consistent.
\end{itemize}
\end{remark}
This shows that the argument described for $\CH$ is far more general and can be used to single out many other set theoretic properties. We can apply it for example to $2^{\aleph_0}>\aleph_2$:
Moore introduced in \cite{MOO06} a $\Pi_2$-sentence $\theta_{\mathrm{Moore}}$ for $\in_{\Delta_0}$ to show the existence of a definable well order of the reals in type $\omega_2$ in models of the bounded proper forcing axiom.\footnote{We use here the definable well-ordering of the reals in type $\omega_2$ existing in models of bounded forcing axioms isolated by Moore, but Thm. \ref{mainthm:2omegageqomega2*} could be proved replacing $\theta_{\mathrm{Moore}}$ with any other coding device which produce the same effects, for example those introduced in \cite{CAIVEL06,TOD02}, or the sentence $\psi_{\mathrm{AC}}$ of Woodin as in \cite[Section 6]{HSTLARSON}.} 

\begin{theorem}  \label{mainthm:2omegageqomega2*}
There is a $\Pi_2$-sentence $\theta_{\mathrm{Moore}}$  for $\in_{\Delta_0}$ such that
the following holds:
\begin{enumerate}
\item \label{mainthm:2omegageqomega2*-1}
$\theta_{\mathrm{Moore}}$ is independent of $S+T_{\Delta_0}$, where $S$ is the $\in$-theory of Thm. \ref{Thm: mainthmforcibility}. 
\item \label{mainthm:2omegageqomega2*-2}
$\ZFC^-_{\Delta_0}+\exists x\,(x\text{ is }\aleph_1)+\theta_{\mathrm{Moore}}$ proves that there exists a well-ordering of\footnote{More precisely: there is a $\ZFC^-_{\Delta_0}$-provably $\Delta_1$-property $\psi(x,y,z)$ such that in any model $\mathcal{M}$ of the mentioned theory there is a parameter $d\in\mathcal{M}$ such that $\psi(x,y,d)$ defines an injection of $\pow{\omega}$ of the model 
with the class of ordinals of size at most $\omega_1$ of the model.} $\pow{\omega}$ in type $\omega_2$.
\item \label{mainthm:2omegageqomega2*-3}
$\ZFC^-_{\Delta_0}+\exists x\,(x\text{ is }\aleph_2)+\theta_{\mathrm{Moore}}$ proves\footnote{Recall that $2^{\aleph_0}=\omega_2$ is the conjunction of the $\Sigma_2$-sentences \emph{$\aleph_2$ exists} and \emph{there is a bijection of $\aleph_2$ onto the reals}, in particular we need the assumption $\exists x\,(x\text{ is }\aleph_2)$ to hold in a fragment of $\ZFC$ which is able to make sense of the sentence $2^{\aleph_0}=\omega_2$. This sentence fails for trivial reasons in a model of \emph{$\aleph_2$ does not exist} even if the model has a definable well order of the continuum in type $\aleph_2$. This is the reason why this item of the Theorem is distinct from the previous one.} that $2^{\aleph_0}=\aleph_2$.
\item \label{mainthm:2omegageqomega2*-4}
For $S$ and $\in_B$ the theory and signature considered in Thm. \ref{Thm: mainthmforcibility},
$\exists x\,(x\text{ is }\aleph_1),\theta_{\mathrm{Moore}}$ are both in $\bool{AMC}(R,B)$ for any $\in$-theory $R$ extending $S$.
\item \label{mainthm:2omegageqomega2*-5}
If $R$ extends $\ZFC$, $A\in\SpecAMC{R}{\in}$  is such that $\in_A$ contains $\in_{\Delta_0}$, $\exists x\, (x\text{ is }\aleph_2)\in \bool{AMC}(R,A)$ and
\[
\theta_{\mathrm{Moore}}+(R+T_{\in,A})^{\in_A}_{\forall\vee\exists}+\ZFC
\]
is consistent, then $2^{\aleph_0}>\aleph_2$ is not in $\bool{AMC}(R,A)$.
\end{enumerate}
\end{theorem}

The two theorems single out $2^{\aleph_0}=\aleph_2$ among all possible solutions of the continuum problem.\footnote{These results bring to light the role of forcing axioms in deciding the value of the continuum; an overview of the proofs of $2^{\aleph_0}=\aleph_2$ given by forcing axioms is given for example in \cite{Moorecont}.} For any $\in$-theory $R$ extending $\ZFC+$\emph{large cardinals} there is at least one $B\in\SpecAMC{R}{\in}$ with $\neg\CH$ (actually with a definable version of $2^{\aleph_0}=\aleph_2$) in $\bool{AMC}(R,B)$, and this occurs even if $R\models\CH$ or $R\models 2^{\aleph_0}>\aleph_2$.
On the other hand for any $\in$-theory $R$ extending $\ZFC$, if $\CH$ is independent of $R$, then $\CH$ is never in $\bool{AMC}(R,A)$ for any $A\in\SpecAMC{R}{\in}$
(if $\in_A$ contains $\in_{\Delta_0}$)
and similarly if $\theta_{\mathrm{Moore}}$ is independent of $R$, $2^{\aleph_0}>\aleph_2$  is never in $\bool{AMC}(R,A)$ for any $A\in\SpecAMC{R}{\in}$ (if $\in_A$ contains $\in_{\Delta_0}$).\footnote{Furthermore the last part of Thm. \ref{Thm: mainthmforcibility} shows that
$\CH$, $2^{\aleph_0}=\aleph_2$, $2^{\aleph_0}>\aleph_2$, $\theta_{\mathrm{Moore}}$ are all 
boolean combination of $\Pi_2$-sentences in the signature $\in_{\Delta_0}$ which cannot be expressed by
boolean combination of $\Pi_1$-sentences for the signature $\in_B\supseteq\in_{\Delta_0}$ in models of $S$ (with $S$ and $B$ as in Thm. \ref{Thm: mainthmforcibility}), 
as their truth value can be changed by forcing.}

\subsection{Model completeness for the theory of $H_{\aleph_2}$ versus bounded forcing axioms and Woodin's axiom $(*)$}

Finally we can show that for certain signatures $\tau$ and a model $(V,\in)$ of $\ZFC+$\emph{large cardinals}, a model completeness result for the $\tau$-theory of $H_{\aleph_2}^V$ is equivalent to assert that $V$ satisfies a strong form of Woodin's axiom $(*)$. The reader interested in the (so far) optimal result is referred to \cite[Thm. 5.6]{viale2021absolute}. Recall the axioms $\BMM$, $\BMM^{++}$, $(*)$ presented in  Section \ref{subsec:exclosedmodst}.
An immediate consequence of the definition of $\in_B$ and of the invariance of the universal $\in_B$-theory of $V$ across forcing extensions is the following result.\footnote{$\in_B$ extends $\in_{\Delta_0}$ with a  a predicate for the non-stationary, and contains a predicate for a (lightface definable in $L(\Ord^\omega)$) set of reals which is Wadge above  any set of reals in $L(\mathbb{R})$.}
\begin{theorem}\label{thm:modcompBMM}
Assume $S$ and $B$ are as in Thm. \ref{Thm: mainthmforcibility} and $(V,\in)$ is a
model of $S$ such that the $\in_B$-theory of $H_{\aleph_2}^V$ is model complete.
Then $(V,\in)$ models $\BMM^{++}$.

Assume furthermore that $\NS_{\omega_1}$ is precipitous in $V$.
Then $(V,\in)$ models $(*)$.
\end{theorem}
Note that $\NS_{\omega_1}$ is a precipitous ideal in $V$, if $V$ models Martin's maximum (of which $\BMM$ is a weakening).
The theorem follows almost immediately from Thm. \ref{Thm: mainthmforcibility} by the observation that if $G$ is $V$-generic for a stationary set preserving forcing and $(V,\in)$ models $S$ \[(H_{\omega_2},\in_B^V)\sqsubseteq(V[G],\in_B^{V[G]}).\]

\section{On maximality and justification}\label{sec:philcons}

In this last section we focus on two philosophical aspects of the present approach.
First we address the following question:
which notion of maximality is displayed by the model companions of $\ZFC +$\emph{large cardinals}?
Then, we study the role of  large cardinals and strong forcing axioms in providing model completeness for the theory of $H_{\aleph_2}$. The role these axioms play in obtaining these model companionship results suggests a second question: is there a form of justification of set-theoretical axioms that emerges from the present approach? In the remainder of this section we address these two questions.

\subsection{Maximality}

Maximality principles play an important role in the study of the axiomatic extensions of $\ZFC$. Maximality conveys the idea that the universe of sets is as full as it could be, while maximality principles are (first order) statements which try to capture this idea axiomatically.\footnote{See \cite{Incurvati2016} for a complete and informative survey on the topic of maximality principles in set theory.} 

There is no unique way to realize maximality in set theory, and its various forms depend on what in fact is maximized. Indeed, the process of maximization can be applied to objects, or possibilities.\footnote{To be faithful to the literature, we should also mention another form: the maximization of domains. This is a strategy inspired by Hilbert's (second order) axiom of completeness that was (almost informally) presented in the following terms: ``The elements (points, straight lines, planes) of geometry form a system of things that, compatibly with the other axioms, can not be extended; i.e. it is not possible to add to the system of points, straight lines, planes another system of things in such a way that in the resulting system all the axioms I-IV, V.1 are satisfied''. Although this axiom was meant to provide the completeness of the real line (and consequently of the real plane), it is standard to interpret its (somewhat vague) formulation as stating, in a second order way,  the maximality of the domain of interpretation of Hilbert's axioms of geometry. Contrary to the other forms of maximality, this approach did not give rise to new well-studied set theoretical principles. However, a first order axiomatic realization of this approach is given by McGee in terms of his Completeness Principle \cite{McGee1992-MCGXPW}. McGee's principle, like Hilbert's axiom, is sufficient to prove categoricity for the class of pure sets (in a context of a set theory with Urelemente).}

\smallskip

\begin{itemize}

\item With respect to objects we find two standard forms of maximization: height maximality (expressing maximality for the lengths of the series of ordinals) and width maximality (expressing the maximality of the power-set operation). The former is successfully exemplified by  large cardinals axioms \cite{Koellner2010-KOEIAL}, while the latter represents a cluster of notions often connected to the justification of forcing axioms, generic absoluteness principles \cite{Bag2005}, and forms of horizontal reflection like the Inner Model Hypothesis \cite{AntosForthcoming-ANTUAE-2}.

\item For what concerns the maximization of possibilities, we find modal principles of forcing like the so-called Maximality Principle \cite{Hamkins2003-HAMASM} (expressing the idea that everything that is possibly necessary is true) or resurrection axioms \cite{Hamkins2014-HAMRAA-4,VIAASP,VIAAUD14} (expressing the idea that if something is true, then it is necessarily possible).  This form of maximality is closely related to forms of generic absoluteness, but it differs from it conceptually. Indeed, the justification of these modal principles often relays on potentialist views of set theory \cite{linnebo_2013,https://doi.org/10.1111/nous.12208} that are not always compatible with the actualist ones that motivate generic absoluteness 
principles.\footnote{Notice that this incompatibility only concerns the philosophical justification of these principles, but does not affect their mathematical coherence with generic absoluteness principles (which can be justified on the basis of width or height maximality).}


\end{itemize}
\smallskip

Now, \emph{what form of maximality is displayed by the notion of model companionship?} Can it be compared to any of the above forms of set-theoretical maximality?

\medskip

A first step toward answering the above questions consists in analyzing whether the maximality provided by models companions concerns primarily the underlying ontology (thus suggesting a form of maximization of objects), or the class of sentences that are considered valid (as it seems the case when considering generic absoluteness principles or modal principles of forcing). In other words, does the maximality provided by model companions concerns primarily the ontology of sets, or the family of set-theoretic truths? 
This question does not have a clear-cut answer. As a matter of fact, we believe that model companionship results for set theory display a multifaceted character related to both forms of maximality.

For what concerns height maximality, this notion is successfully realized by large cardinal axioms. While the model companionship results for set theory do not seem to enforce this form of maximality, they can nonetheless be used to justify it (we argue on this point in the next section).
On the other hand, for what concerns width maximality, notice that bounded forcing axioms (and similar principles meant to express forms of horizontal maximality) can be reformulated in terms of generic absoluteness properties for $H_{\omega_2}$. Many of the results we presented in the previous section outlines that such generic absoluteness results come in pairs with the model completeness for the theory of $H_{\aleph_2}$ with respect to natural signatures for this structure (e.g. Thm. \ref{thm:modcompBMM} or \cite[Thm. 5.6]{viale2021absolute}). Note that the standard formulation of bounded forcing axioms, such as $\BMM,\BMM^{++},(*)$, yield forms of absoluteness for the theory of $H_{\aleph_2}$, which are (apparently) weaker than the model completeness result given by Thm. \ref{Thm: mainthmforcibility}:
they only  assert that
(for the appropriate fragment $\tau$ of $\in_B$)
$H_{\aleph_2}$ is $\Sigma_1$-elementary for $\tau$ only with respect to certain forcing extensions of $V$
(which are $\tau$-models of the universal $\tau$-theory of $V$ by the second part of Thm. \ref{Thm: mainthmforcibility}). On the contrary, model completeness of the $\in_B$-theory of $H_{\aleph_2}$ yields $\Sigma_1$-elementarity for the much larger class of \emph{all} $\in_B$-superstructures which satisfy the same universal $\in_B$-theory.


For what concerns the modal principles of forcing, we can argue as for width maximality. As a matter of fact, the Maximality Principle and resurrection axioms, as exposed in \cite{VIAASP,Hamkins2003-HAMASM,Hamkins2014-HAMRAA-4,VIAAUD14} are, by all means, principles of  generic absoluteness for the theory of $H_{\aleph_2}$ (at least in the presence of mild forms of bounded forcing axioms). 
Indeed, for many of these principles their natural formulation asserts that the theory of $H_{\aleph_2}$ is given by those $\Pi_2$-sentences $\varphi$ (for an appropriate fragment of $\in_B$) for which $\varphi\to\square\Diamond\varphi$ holds in the (Kripke frame given by the) generic multiverse composed of initial segments of forcing extensions of $V$. Therefore, they describe the theory of $H_{\aleph_2}$ with respect to the generic multiverse much like  infinitely generic structures (introduced by Robinson in connection with the notion of infinite forcing) describe the theory of the existentially closed models for a first order theory. Leveraging on the model completeness results for $\in_B$ we can actually identify the theory of $H_{\aleph_2}$ in models of strong forcing axioms with that given by the infinitely generic structures, in signature $\in_B$, for the class of models of $\ZFC+$\emph{large cardinals}. Hence the absoluteness displayed by the model companions of set theory, relative to $\in_B$, largely encompasses that of the modal principles of forcing, since the latter are infinitely generic (at least on the face of their definition) only with respect to the generic multiverse and not with respect to the elementary class of the $\in_B$-theory of $\ZFC_B+$\emph{large cardinals}.\footnote{See \cite{Venturi2020-VENIFA} for the details on this comparison between infinite forcing and modal principles of forcing.}

Finally, what is the connection between the maximization given by model companionship results and the (somewhat vague) form of completeness that is encoded in the principle that we called Hilbertian Completeness?  Although our aim, similar to Hilbert's, is  to implement the idea that the universe $(V, \in)$ contains \emph{all possible}  sets, our formal way to make this idea concrete differs significantly from Hilbert's approach. Instead of focusing directly on the objects, we maximize existential sentences using our principle of Informal Existential Completeness \ref{Descr:IWF}. Thus, maximality of the objects receives in our approach a syntactic twist, being formulated in terms of maximality of solutions to simple problems (according to a syntactic measure of simplicity given by the selection of a specific first order signature). 

Summing up, it seems that the form of maximality provided by (absolute) model companionship does not precisely fit with any of the standard forms of set-theoretical maximality. Nonetheless there are very strong similarities with generic absoluteness principles and certain types of forcing axioms\footnote{To the extent that a strong form of Woodin's axiom $(*)$ can be equivalently formulated as the assertion that for a certain signature $\tau$ the $\tau$-theory of $H_{\aleph_2}$ is the AMC of the $\tau$-theory of $V$.}, and conceptual similarities with Hilbert's idea of completeness. We see the notion of model companionship as a useful tool towards a mathematical formalization of the somewhat vague ideas of width maximality and Hilbertian completeness. Our opinion is that the model companions of set theory yield a form of maximality that we may call \emph{algebraic maximality}. The inspiration clearly comes from algebraically closed fields. We consider algebraic maximization of a theory $T$ as a two steps process: first, one determines which are the basic, simple, concepts of the theory of interest; then, once the language is extended  to a signature $\tau$ with symbols for these concepts, one analyzes the $T$-ec $\tau$-models. $T$ is \emph{algebraically maximal for $\tau$} in case the $\tau$-theory of the
$T$-ec $\tau$-models is model complete and also when its $T$-ec $\tau$-models can be seen as ``parts of'' the models of $T$. This can be formally implemented as follows:
\begin{definition}
Let $T$ be a $\tau$-theory and $A$ be a subset of $\bool{Form}_\tau\times 2$.
$T$ is \emph{algebraically maximal} relative to $A$, if $A\in\SpecAMC{T}{\tau}$, and
there is a $\tau$-formula $\psi_A(x)$ such that 
for any $\mathcal{M}=(M,\tau^M)$ which models $T$:\footnote{$\psi_A(\mathcal{M})$ is the extension of $\psi_A(x)$ in $\mathcal{M}$, and $\tau_A^{M}$ is the unique interpretation of the predicates of $\tau_A$ in $M$ which can make $(M,\tau_A^M)$ an expansion of
$\mathcal{M}$ which is a $\tau_A$-model of $T_{\tau,A}$.}
\begin{itemize} 
\item
$(\psi_A(\mathcal{M}),\tau_A^{M})$  realizes $\bool{AMC}(T,A)$,
\item
$(\psi_A(\mathcal{M}),\tau_A^{M})$ is a $\tau_A$-substructure of $(M,\tau_A^M)$.
\end{itemize}
\end{definition}

In terms of  algebraic maximality we briefly review the  results of Section \ref{sec:mainresults}. First of all notice that
Thm. \ref{Thm:AMCsettheory+Repl} states that for any extension $T$ of $\ZFC$, algebraic maximality can be obtained for $T$ by sets $A_\kappa$ such that the AMC of $T$ in signautre $\in_{A_\kappa}$ describes the theory of $H_{\kappa^+}$ in models of $T$, as $\kappa$ ranges among the $T$-definable cardinals.\footnote{Notice that we do not assert that the set theoretic concepts made elementary by the different $A_\kappa$'s can be considered simple according to any ``reasonable'' simplicity criterion for $\kappa$. We give separately below arguments showing that for $\kappa=\aleph_0$ or $\kappa=\aleph_1$, $A_\kappa$ can be chosen to be ``simple'' according to our simplicity criteria. As of now, we do not have anything to say about the ``simplicity'' of $\in_{A_\kappa}$ for any $\kappa>\aleph_1$.}
Moreover, the first part of Thm. \ref{Thm: mainthmforcibility} states that algebraic maximality, relative to $B$, holds for the theory $\ZFC+$\emph{large cardinals}$+$\emph{forcing axioms}, while its second part states that $\in_B$ satisfies the simplicity criteria set forth in \ref{Descr:IWF} of Section \ref{sec:AMCspecST}.
Finally, we have that Thm. \ref{thm:modcompBMM} (and its optimal strengthening given by \cite[Thm. 5.6]{viale2021absolute})
shows that algebraic maximality for the $\in_B$ of  Thm. \ref{Thm: mainthmforcibility} entails (or is equivalent to) certain strong forcing axioms.

In conclusion we can argue that Thm. \ref{Thm: mainthmforcibility} and 
Thm. \ref{thm:modcompBMM} show that 
algebraic maximality for set theory,  relative to $\in_B$, describes the \emph{correct} axioms for the theory of $H_{\aleph_2}$,
as algebraic maximality for $\in_B$ realizes \ref{Descr:IWF} for this fragment of the universe. This naturally lead us to discuss the second topic of this section: justification.

\subsection{Justification}

Let us now turn to the topic of justification. When it comes to justifying new set-theoretical principles, the proposed reasons normally fall into two main categories: \emph{intrinsic} or \emph{extrinsic}. These forms of justification have been developed within the  so-called G\"odel's program\footnote{Presented by G\"odel in his seminal paper \cite{God47a}.}: 
a step by step extension of $\ZFC$, aimed to coherently complete our picture of the universe of sets. In this sense, the present approach is clearly in the wake of G\"odel's program. Indeed, the axiomatization of the model companions of set theory should, eventually, provide a description of the theory of $V$ as the stratification of the theories of the various $H_{\kappa^+}$. 
Furthermore, the essential role that large cardinals and strong forms of forcing axioms play in the individuation of these model companionship results suggests a justification of the truth of these axioms in view of the nice model-theoretic properties for the theories of $H_{\aleph_1}$ or $H_{\aleph_2}$ they entail.

In order to better argue for this form of justification of large cardinals and forcing axioms, let us briefly revise what we mean by intrinsic and extrinsic reasons. 

\medskip

\begin{itemize}

\item  A form of justification is considered intrinsic when it is based on intrinsic features of the concept of set or, derivatively (given the foundational role of set theory), on notions that are key to the whole edifice of mathematics and logic. In this sense, a new axiom is justified when it captures an essential aspect of the concept of set or when it formalizes notions that are fundamental for mathematics and logic. This form of justification is utterly conceptual, since it rests on the theoretical priority of the notion of set over its formalization. Consequently, an extension of $\ZFC$ is well justified, when it faithfully represents the correct concept of set (e.g. the iterative conception \cite{Boolos1971-BOOTIC} or the quasi-combinatorial one \cite{Bernays1935-BERSLP-9}). 

\item A form of justification is considered extrinsic when we are forced to accept an axiom by the abundance of its desirable consequences; even in the absence of intrinsic reasons. In this case an axiom is justified even if it does not seem \emph{prima facie} to capture any relevant aspect of the  concept of set. This second form of justification is clearly meant to overcome the limits of intrinsic reasons and to account for an experimental methodology in the context of the foundations of mathematics. Extrinsic justifications are pragmatic in nature, since they rely on the fruitfulness of an axiom in solving new problems, shortening proofs of already known results, and unifying substantial bodies of theory. The form of justification put forward by extrinsic reasons is akin to an inference to the best explanation within mathematics. Therefore, an extrinsically justified axiom is judged by its consequence and not by its meaning. 

\end{itemize}

\medskip

Intrinsic and extrinsic reasons have been extensively studied in the philosophy of set theory \cite{Maddy1997-MADNIM,Maddy2007-MADSPA-2} and they have been widely applied, in recent debates, for the justifications of competing programs  \cite{10.1093/philmat/nkj009,Arrigoni2013-ARRTHP}. These two forms of justification have also been criticized for their opacity in offering clear criteria of application and for the lack of demarcation between intrinsic and extrinsic reasons \cite{Barton2020-BAROFO}. 
Moreover, as it happened in the debate on the analytic-synthetic distinction, 
there is also no shortage of contributions that reject the problem of justification at its very base. Following a naturalistic account of mathematics, authors like Hamkins (or before him Cohen \cite{10.1216/rmjm/1181070010}) propose to dismiss the problem of justification, together with the problem of independence, by defending the liberty of mathematicians to study different universes of set theory and by declaring the study of the variation of truth values among the models of $\ZFC$ to be all that mathematically matters for the study of independence \cite{Hamkins}. 

\medskip

Now, \emph{in which sense the nice model-theoretic properties of model companions are able to justify large cardinals and strong forcing axioms}? Is the role of the latter set-theoretical principles in the individuation of model companions of set theory able to provide an intrinsic or an extrinsic justification for them? 

\medskip

A first complication that we face in addressing these questions is that the present approach does not only deal with axiomatic extensions of $\ZFC$, but also with linguistic ones: the introduction of new signatures for set theory. This peculiar aspect is responsible for a justification of large cardinals and strong forcing axioms that  contains elements of both intrinsic and extrinsic arguments. Indeed, on the one hand strong forcing axioms provide a maximization of the $\Pi_2$-statements realized over $H_{\aleph_2}$ and, thus,  have tremendously abundant consequences on third order arithmetic (and on this ground they can be extrinsically justified). On the other hand, because of the possibility to provide (generic) absoluteness results, large cardinal axioms determine the ``simple'' concepts that need to be included in the new signatures for set theory.\footnote{The difficulty in choosing whether forcing axioms can be justified on intrinsic or extrinsic grounds in not new. See \cite{Venturi2019} for a discussion of the connection between forcing axioms and their intrinsic justification in terms of the clarification of the notion of arbitrary set.}  Therefore, large cardinals and strong forcing axioms help us understanding what are the basic concepts on which third order arithmetic should be based (and for this reason they can be intrinsically justified).\footnote{To be more precise, the justification of large cardinals that we are proposing only concerns those axioms of infinity that appear in the model companionship results for third order arithmetic (i.e. large cardinals in the range of supercompactness). In this sense this claim is weaker than it seems. The situation is similar to what happens when we justify Woodin cardinals on the basis of the generic absoluteness results for $L(\mathbb{R})$ (we thank one of the referees for pointing out this analogy). Our argument for justifying large cardinals on the basis of the model companionship results, however, gives an explicit and mathematically sound formulation of the intrinsic justification of large cardinals rooted on the -somewhat vague- assertion that these axioms render ``simple'' many ``natural'' set theoretic concepts.}

Another aspect of the present approach (that places the justification of new axioms somewhat outside the standard practice) is the focus on models, instead of sets (like intrinsic reasons), or problems (like extrinsic ones). As a matter of fact,  large cardinals and strong forcing axioms play an important role in the construction of model companions of set theory, by detecting the natural signatures for which one has nice model-theoretic properties for (the models of) set theory. Consequently, the explicit link these axioms establish between the choice of a signature and model companionship results for set theory can also be seen as an indirect effect of  these axioms in capturing certain natural features of the universe of sets. In this sense,  the present approach makes evident the role that certain additional axioms of set theory  have in clarifying the model theory of set theory, but also, and derivatively, the role that model-theoretic properties have in selecting what could be the natural theories for initial segments of the universe of sets. It is with respect to this double role (selecting simple concepts and providing natural theories) that our approach rooted in model companionship argues for the justifications of new axioms of set theory. 

 Because of the limits of the standard intrinsic-extrinsic dichotomy to capture the form of justification of this model-theoretic approach, we may call this new form of justification  \emph{algebraic}, since it rests on the model-theoretic properties of algebraically maximal theories. It is a form of justification that has both intrinsic and extrinsic elements and that heavily relies on the model-theoretic form of maximality that we described in the previous section. 

In order to exemplify this algebraic form of justification, let us see how the present approach advances G\"odel's program in providing well-justified extensions of $\ZFC$. 
This particular implementation of G\"odel's program is devised to provide new well-justified axioms of larger and larger initial segments of the universe $V$, concretely, by providing axioms able to make the theories of $H_{\kappa^+}$ (as $\kappa$ ranges among the infinite cardinal) the model companions of set theory for the ``correct'' signature aimed at describing the elementary properties of $\pow{\kappa}$. In this way the theory of $V$ is gradually completed by adding new principles whose justification rest on the nice model-theoretic properties they provide for the theories of the different $H_{\kappa^+}$. Let us see how this justification process can be implemented.

We start, of course, from $\ZFC$ and notice that  this theory is algebraically maximal relative to $\in_A$, which is the signature obtained by adding to $\in_{\Delta_0}$ predicates for the lightface definable subsets of $\pow{\omega}$ in the Chang model $L(\Ord^\omega)$.\footnote{Equivalently one can obtain $\in_A$ from the signature $\in_B$ of Thm. \ref{Thm: mainthmforcibility} removing the constant symbol for $\omega_1$ and the predicate symbol for the non-stationary ideal on $\omega_1$.} In this case its existentially complete models describe the theory of $H_{\omega_1}$, and -because $H_{\omega_1}$ is a definable sub-model of $V$- we have that all condition of algebraic maximality are fulfilled. 
By  \cite[Thm. 5.2]{VIAVENMODCOMP}, this theory is the AMC of any extension of $\ZFC$ for $\in_A$.

However the set theoretic concepts expressed by atomic $\in_A$-formulae or by the universal $\in_A$-theory of $\ZFC+T_{\in,A}$ are not simple according to any of our criteria (put forward in  \ref{cond:simplstprop}  of Section \ref{sec:AMCspecST}).  For example $\forall x\,(x\subseteq\omega\to x\in L)$ is a $\Pi_1$-sentence for $\in_A$ which holds in  the $\ZFC$ model given by the constructible universe $L$ and fails in $L[x]$ where $x$ is a Cohen real over $L$. On the other hand, assuming as our base set theory 
$\ZFC+$\emph{there are class many Woodin cardinals which are a limit of Woodin cardinals}, the 
$\in_A$-axiomatization of this set theory satisfies our simplicity criteria.
Henceforth, we should accept large cardinal axioms
-among other reasons- also because they render simple (according to our criteria) the notions codified by $\in_A$. Consequently, on the base of this argument one can claim that the AMC of set theory enriched with large cardinal axioms, in the signature $\in_A$ describes the correct theory of $H_{\aleph_1}$, and thus realizes the first step of our proposed way to implement G\"odel's program. 

A similar argument can be applied to the theory of $H_{\aleph_2}$, by appealing to Thm. \ref{Thm: mainthmforcibility}. In this case the extension of $\ZFC$ by large cardinal axioms is able to make the signature $\in_B$, in which Thm. \ref{Thm: mainthmforcibility} is expressed, simple according to our criteria of simplicity. Furthermore the nice model-theoretic properties of the theory of  $H_{\aleph_2}$ described by the AMC of $\ZFC+$\emph{large cardinals} relative to $\in_B$ renders $\ZFC+$\emph{large cardinal axioms$+$strong forms of bounded forcing axioms} algebraically maximal with respect to the signature $\in_B$. We can thus justify forcing axioms and large cardinals as they provide (by means of algebraic maximality relative to $\in_B$) the correct picture of the theory of $H_{\aleph_2}$. 

Now, starting with base theory $\ZFC+$\emph{large cardinals}$+$\emph{strong forms of bounded forcing axioms}, we should search for criteria to detect which definable subsets of $\pow{\omega_2}$ describe simple set theoretic notions, and search for the AMC of $\ZFC+$\emph{large cardinals}$+$\emph{strong forcing axioms} relative to the signature expanded by these concepts.
If and when this task is achieved,\footnote{To be honest, we currently do not have any clue on how this task could be performed.} we can argue that the corresponding AMC provides the correct description of $H_{\aleph_3}$ and that the relevant axioms should be added to the axioms of set theory.

In a nutshell this is the general strategy: to climb the cardinal hierarchy in order to implement algebraic maximality for set theory with respect to signatures that describe the properties of $\pow{\kappa}$ for larger and larger cardinals $\kappa$. This will come in pairs with set-theoretic axioms that fulfill (possibly weaker and weaker) criteria of simplicity for these signatures\footnote{For example $\in_A$ does not satisfy the simplicity criterion \ref{cond:simplstprop(ii)}  of page \pageref{cond:simplstprop(ii)} relative to $\ZFC$, while it does relative to $\ZFC+$\emph{large cardinals}; henceforth as our theories increase in logical strength, more and more concepts could become simple according to the criteria set forth in page \pageref{cond:simplstprop(ii)}. Note that $\in_{\Delta_0}$ would have been recognized as a natural signature codifying by atomic formulae simple set-theoretic concepts even before the discovery of forcing and the introduction of large cardinal axioms; on the other hand, $\in_A$ would have not been considered natural before the generic absoluteness results for second order arithmetic and the clear picture large cardinals give of the theory of universally Baire sets. As it occurred for forcing invariance, it might be the case that, as we accept as valid set-theoretic truths certain set-theoretic principles, other criteria of simplicity may emerge.} (which are algebraically maximal for more and more of their initial segments). This strategy has been successful so far with respect to $H_{\aleph_1},H_{\aleph_2}$. Thm. \ref{Thm:AMCsettheory+Repl} shows that no obstruction to its implementation for other cardinals is in sight.\footnote{When the concepts making the theory of $H_{\kappa^+}$ an AMC of set theory are simple, one should expect this AMC to be ``empirically'' complete, i.e. to decide most ``relevant'' set-theoretic problems whose solution is correctly computed in $H_{\kappa^+}$. Again this ``empirical completeness'' holds for the AMCs describing the theories of $H_{\aleph_1},H_{\aleph_2}$.}

We conclude by noting that there is another sense in which these model companionship results, the justification they offer to new axioms,  and the solution they yield for the continuum problem are utterly Hilbertian. This is the autonomy they provide for mathematics with respect to its foundations. As a matter of fact, the use of the notion of model companionship to fix the theory of second and third order arithmetic (and to enforce algebraic maximality)  does not introduce any element that is foreign to mathematical practice. The sought solution to  the continuum problem is a solution obtained with formal tools and that is justified on a purely mathematical ground. 
In this sense the results presented in this paper seem to realize in full the possibility of an autonomous foundation for mathematics, as the one sought by Hilbert \cite{Franks2009-FRATAO-19}. 


\bibliographystyle{plain}
	\bibliography{Biblio}

\end{document}